\documentclass{IEEEtran}
\usepackage{cite}
\usepackage{amsmath,amssymb,amsfonts}
\usepackage{algorithmic,color}
\usepackage{epstopdf}
\usepackage{graphicx}
\usepackage{textcomp}
\usepackage{mathrsfs}
\def\BibTeX{{\rm B\kern-.05em{\sc i\kern-.025em b}\kern-.08em
    T\kern-.1667em\lower.7ex\hbox{E}\kern-.125emX}}

\newtheorem{theorem}{Theorem}
\newtheorem{corollary}{Corollary}
\newtheorem{remark}{Remark}
\newtheorem{lemma}{Lemma}
\newtheorem{assumption}{Assumption}

\usepackage{algorithm}
\usepackage{algorithmic}
\allowdisplaybreaks
\begin{document}

\title{Distributed Solving of Linear Quadratic Optimal Controller with Terminal State Constraint}

\author{Wenjing Yang, Zhaorong Zhang, Juanjuan~Xu
\thanks{This work was supported by the National Natural Science Foundation of China under Grants 62573262, 62503289 and the Natural Science Foundation of Shandong Province under Grants ZR2021JQ24. \textit{(Corresponding author: Juanjuan~Xu.)}}
\thanks{Wenjing Yang is with School of Control Science and Engineering, Shandong University, Jinan, Shandong, P.R. China 250061.
        {\tt\small yangwenjing1024@163.com}}%
\thanks{Zhaorong Zhang is with School of Computer Science and Technology, Shandong University, Qingdao, Shandong, P.R. China 266237.
        {\tt\small zhangzr@sdu.edu.cn}}%
 \thanks{Juanjuan~Xu is with School of Control Science and Engineering, Shandong University, Jinan, Shandong, P.R. China 250061.
        {\tt\small juanjuanxu@sdu.edu.cn}}%
}

\maketitle

\begin{abstract}
This paper is concerned with the linear quadratic (LQ) optimal control of continuous-time system with terminal state constraint. In particular, multiple agents exist in the system which can only access partial information of the matrix parameters. This makes the classical solving method based on Riccati equation with global information suffering. The main contribution is to present a distributed algorithm to derive the optimal controller which is consisting of the distributed iterations for the Riccati equation, a backward differential equation driven by the optimal Lagrange multiplier and the optimal state. Furthermore, the proposed distributed iteration method is extended to solve the consensus control problem for heterogeneous multi-agent systems, achieving the globally optimal performance of the system. The effectiveness of the proposed algorithm is verified by two numerical examples, where the performance index under the proposed distributed controller is smaller than that under the commonly used consensus control.
\end{abstract}


%
\IEEEpeerreviewmaketitle

\section{Introduction}

Linear quadratic (LQ) optimal control is one of the core fundamental problem in modern control theory which has wide applications in plenty of fields including power system scheduling \cite{1}, economics \cite{2} and autonomous driving \cite{3}. A typical example is the trajectory tracking control problem which can be solved by designing a LQ-based tracking controller \cite{4}. For the classic LQ optimal control problem, one of the most powerful tools is Riccati equation \cite{5}, which can be used to design the optimal feedback gain. For the irregular LQ optimal control, a hierarchical optimization method has been proposed in \cite{6} where two Riccati equations are given to design the optimal controller.

Considering the wide existence of terminal state constraints in lots of practical problems such as the precise terminal guidance in the hypersonic vehicle re-entry tracking control \cite{7}, it is necessary to investigate the LQ optimal control problems with terminal state constraints.
\cite{8} studied the nonlinear systems with terminal constraints by using a robust receding horizon controller.
\cite{9} considered a time-varying set of terminal constraints by using robust constrained model predictive control algorithm.
\cite{10} studied LQ optimal control with fixed terminal states and integral quadratic constraints. For the cases with stochastic uncertainty, \cite{11} studied a general class of stochastic LQ control problems with terminal states constrained to convex sets. \cite{12} investigated a stochastic LQ optimal control with a fixed terminal state and linear manifold constraint on the initial state. \cite{13} obtained the explicit optimal solution for discrete-time stochastic LQ optimal control with initial and terminal state constraints.

Note that the design of the optimal controllers in the above mentioned results relies on the global information of the matrix parameters in the systems. This loses efficacy when there exist sensor networks composed of a large number of sensor and actuator nodes/agents with communication capabilities among them \cite{14}, which have received a great deal of attention in the past decade due mainly to its current and future envisioned applications including mobile vehicles, smart grids, environmental monitoring and industrial process control \cite{15,16,17,18}. In this case, each agent can only access partial information of the system and make decisions by using their own and/or neighbors' information, which has given rise to decentralized and distributed control problems.

Decentralized control refers to situations where different control posts access different information and make decisions by using own information. \cite{19} provided the optimal solutions to decentralized linear-quadratic-Gaussian (LQG) control problems with one-step delayed information sharing pattern by using dynamic programming. \cite{20} derived the optimal linear solution of the decentralized LQG control with $d$-step delayed information sharing pattern. The partially nested decentralized LQG problems with controllers connected over a graph were studied in \cite{21}. \cite{22} studied teams with only control sharing, \cite{23} studied the broadcast information structure, and \cite{24} investigated the systems with common and private observations.

Distributed control refers to situations where different agents access different information, share their information through network topology, and make decisions by using own and neighbors' information. One of the most popular problems is consensus control which has been well studied \cite{25, 26}. For the optimal control problems, \cite{27} studied the distributed optimal control for discrete-time systems. \cite{28} discussed the distributed consensus optimal control strategy for multi-agent systems with a general directed graph. \cite{29} studied the design of suboptimal distributed control for the identical decoupled dynamical systems. \cite{30} considered the online distributed convex optimization problem with local dynamic loss function sums under unbalanced graphs. A distributed LQG optimal control problem was investigated in \cite{31}, where each controller has its own input matrix.

Although significant progress has been made for the decentralized and distributed optimal control problem, the following issues still exist:
\begin{itemize}
  \item The matrix parameters in the systems are assumed to be completely known. While the case that each agent can only access partial information of the matrix parameters remains to be solved.
  \item Strong conditions are imposed on the topology and/or the weighting matrices in the cost function to guarantee the distributed structure of the optimal solution, apart from the connectivity. While when these conditions are not established, how to derive the optimal solution in a distributed way remain to be studied.
  \item No terminal state has been imposed in the distributed optimal control problem. While the existence of the terminal constraint additionally introduces an optimal Lagrange multiplier and a backward differential equation driven by the multiplier, whose distributed solving needs further study.
\end{itemize}

In this paper, we will study the LQ optimal control problem with terminal state constraint where there exist multiple agents and each agent accesses partial information of the matrix parameters. The main contribution is to present a distributed algorithm to derive the optimal controller. Specifically, the distributed algorithm is consisting of the distributed iterations for the Riccati equation, a backward differential equation driven by the optimal Lagrange multiplier and the optimal state. Furthermore, we solve the optimal consensus control problem using the proposed distributed iteration method. Numerical examples demonstrate the effectiveness of the proposed algorithms, where the performance index under the proposed distributed controller is smaller than that under the commonly used consensus control.

The remainder of the paper is organized as follows. Section II  describes the studied problem. The explicit and the numerical solutions in centralized way are provided in Section III. Section IV gives the distributed algorithm for the optimal controller. Section V gives the distributed solution to the optimal consensus control problem. Numerical examples are given in Section VI. The conclusion is outlined in the last Section. Detailed proofs of the results are given in Appendix.

The following notations will be used throughout this paper: $R^n$
denotes the set of $n$-dimensional vectors;
$x'$ denotes the transpose of $x$. $I$ is identity matrices with compatible dimensions.


\section{Problem formulation}

Consider the following continuous-time system:
\begin{eqnarray}
\dot{x}(t)=Ax(t)+Bu(t),\label{ds1}
\end{eqnarray}
where $x(t)\in R^n$ is the state, $u(t)\in R^m$ is the control input.
$A$ and $B$ are constant matrices with compatible dimensions.
$t\in [0, T]$ is the continuous-time variable.
The initial state is prescribed as $x(0)=x_0$. The terminal state is predetermined as $x(T)=x_T$.

The cost function is given by
\begin{eqnarray}
J=&~\int_{0}^T[x'(t)Qx(t)+u'(t)Ru(t)]dt,\label{ds2}
\end{eqnarray}
where $Q$ is a positive semi-definite matrix with compatible dimension, and $R$ is a positive-definite matrix with compatible dimension.

In this paper, we study a novel case with a sensor network in the system that consists of $N$ agents or nodes, and the communication network between the agents is represented by the undirected topological graph $\mathscr{G}=\{\mathscr{V},\mathscr{E}, \mathscr{A}\}$. In particular, each agent $i, i=1, 2, \ldots, N$ can only communicate with its neighbors $N_{i}$. The graph is assumed to be connected. For more details about the graph refer to \cite{28}.

In this case, the matrix parameters $A, B$ can not be all available to the controller design of agent $i$. That is, different partial information is accessible to different agent which makes the design of the optimal controller much more complicated. In details, the following assumptions are given to show the detailed information available to agent $i$.

\begin{assumption}
There exist matrices $A_i, B_i, i=1, \ldots, N$ such that $A=\frac{1}{N}\sum_{i=1}^NA_i$, $B=\frac{1}{N}\sum_{i=1}^NB_i$.
In particular, $A_i$ and $B_i$ are accessible for agent $i, i=1,\ldots, N$.
\end{assumption}

\begin{assumption}
There exist matrices $Q_i, i=1, \ldots, N$ such that $Q=\frac{1}{N}\sum_{i=1}^NQ_i$, and $Q_i$ is accessible for agent $i, i=1,\ldots, N$. The matrix $R$ is accessible to each agent.
\end{assumption}

\begin{assumption}
There exist $x_{i,0}, x_{i,T}, i=1, \ldots, N$ such that $x(0)=x_{0}=\frac{1}{N}\sum_{i=1}^Nx_{i,0}$ and $x(T)=x_{T}=\frac{1}{N}\sum_{i=1}^Nx_{i,T}$.
In particular, $x_{i,0}, x_{i,T}$ accessible to agent $i, i=1,\ldots, N$.
\end{assumption}

\begin{assumption}
There exists matrices $M_i, i=1, \ldots, N$ such that $BR^{-1} B'=\frac{1}{N}\sum_{i=1}^NM_iM_i'$, where $M_i$ is accessible for agent $i, i=1,\ldots, N$.
\end{assumption}

Based on the above description, the studied problem is concluded as follows.

\textbf{Problem 1.} The aim of the paper is to design a distributed algorithm for the optimal controller which minimizes the cost function (\ref{ds2}) subject to system (\ref{ds1}) and ensures the state to satisfy $x(0)=x_0$ and $x(T)=x_T$ for each agent $i$ under Assumptions 1-4.

\begin{remark}
It is noted that Assumptions 1-4 are meaningful and have practical applications in network systems \cite{31} and multi-agent systems, where the application of multi-agent systems is studied in Section V. In addition, we also note that although the system matrix was assumed to have a decomposable form under Assumptions 1-4, solving Problem 1 cannot be simply decomposed into the sum of $N$ subproblems, for example, in distributed optimal control problem where Assumptions 1-4 hold naturally, the derivation of the optimal solution is NP-hard \cite{29}.
\end{remark}

\section{Preliminaries}
\subsection{Explicit solution by using all the information}

In this part, we present the optimal solution to the Problem 1 by assuming all the information of matrix parameters are available for the agents.
Noting that the terminal state is predetermined, thus it is necessary to ensure that the system can be driven from the initial state to the terminal state. To this end, we firstly present the definition of the reachability of system (1) as below.

\textit{Definition 1:} A state $x_T$ is said to be reachable at time $T$ from the initial state $x_0$ of system (\ref{ds1}) if there exists a controller $u(t)$, $t\in [0,T]$ such that the solution $x(t)$ of the system (1) satisfies $x(T)=x_{T}$.

Next, we introduce the Riccati differential equation:
\begin{eqnarray}
\dot{P}(t)+A'P(t)+P(t)A+Q-P(t)BR^{-1}B'P(t)=0,\label{ds3}
\end{eqnarray}
with terminal value $P(T)=0$ and the following matrix differential equation:
\begin{eqnarray}
\dot{\Phi}(t)=&~[A-BR^{-1}B'P(t)]\Phi(t),\label{ds4}
\end{eqnarray}
with initial value $\Phi(0)=I$. Furthermore, we denote $\Phi(t,s)=\Phi(t)\Phi^{-1}(s)$ and $\Psi(t)={\Phi'(t)}^{-1}$.

The optimal solution is then given as follows.

\begin{theorem}\label{t1} If the terminal state $x_T$ is reachable from the initial state $x_0$ of system (\ref{ds1}), then there exists a solution $\lambda^*$ to the following equation:
\begin{eqnarray}
\Big[\int_0^T\Phi(T,s)BR^{-1}B'\Psi(s,T)ds\Big]\lambda
=\Phi(T,0)x_0-x_T.\label{ds5}
\end{eqnarray}
Then, the optimal controller of minimizing (\ref{ds2}) subject to (\ref{ds1}) with $x(0)=x_0$ and $x(T)=x_T$ is given by:
\begin{eqnarray}
u^{*}(t)=-R^{-1}B'P(t)x^{*}(t)-R^{-1}B'\beta(t),\label{ds6}
\end{eqnarray}
where $\beta(t)$ obeys the following equation:
\begin{eqnarray}
\dot{\beta}(t)=-[A'-P(t)BR^{-1}B']\beta(t),~\beta(T)=\lambda^*,\label{ds8}
\end{eqnarray}
and $x^{*}(t)$ obeys the following equation:
\begin{eqnarray}
\dot{x}^{*}(t)&=&[A-BR^{-1}B'P(t)]x^{*}(t)-BR^{-1}B'\beta(t),\label{ds9}
\end{eqnarray}
with $x^{*}(0)=x_{0}$.
\end{theorem}

\textit{Proof:} The proof follows directly from Theorem 1 in \cite{33} and is thus omitted here.  \hfill $\blacksquare$

\subsection{Numerical solution by using all the information}

From Theorem \ref{t1}, the design of optimal controller (\ref{ds6}) depends on the solution $P(t)$ of the Riccati differential equation (\ref{ds3}), which is nonlinear and typically challenging to solve analytically. To address this, we propose the following iteration:
\begin{eqnarray}
\dot{P}^{n}(t)&=-Z^{n}(t)'P^{n}(t)-P^{n}(t)Z^{n}(t)-V^{n}(t), \label{1}
\end{eqnarray}
where $t\in [0, T]$ and
\begin{eqnarray}
Z^{n}(t)&=&A-BR^{-1}B'P^{n-1}(t),\label{7}\\
V^{n}(t)&=&Q+P^{n-1}(t)BR^{-1}B'P^{n-1}(t),\label{8}
\end{eqnarray}
with $P^{n}(T)=0$, $P^{0}(t)=0$, and give its solution by
\begin{eqnarray}
P^{n}(t)=\Psi^{n}(t)\big[\int_t^T\Phi^{n}(s)'
V^{n}(s)\Phi^{n}(s)ds\big]\Psi^{n}(t)',\label{24}
\end{eqnarray}
where $\Phi^{n}(t)$ and $\Psi^{n}(t)$ satisfy the following iteration respectively:
\begin{eqnarray}
\dot{\Phi}^{n}(t)&=&Z^{n}(t)\Phi^{n}(t), \label{2}\\
\dot{\Psi}^{n}(t)&=&-Z^{n}(t)'\Psi^{n}(t), \label{3}
\end{eqnarray}
with $\Phi^{n}(0)=I$ and $\Psi^{n}(0)=I$.

It is then easy to verify the convergence of the iteration (\ref{1}) as given below.

\begin{lemma}\label{lem1}
It holds that

(1) $P^{n}(t)\geq P^{n+1}(t)$;

(2)$\displaystyle\lim_{n\rightarrow \infty} P^{n}(t)= P(t)$,
where $P(t)$ satisfies (\ref{ds3}).
\end{lemma}

 \emph{Proof:} The proof can be directly obtained from Theorem 7.2 in \cite{34} and is thus omitted here.  \hfill $\blacksquare$

Furthermore, it is immediate to obtain the following result.
\begin{corollary}\label{cor1} It holds that
\begin{eqnarray}
\lim_{n\rightarrow \infty}Z^{n}(t)&=&Z(t), \\
\lim_{n\rightarrow \infty}\Phi^{n}(t)&=&\Phi(t),\\
\lim_{n\rightarrow \infty}\Psi^{n}(t)&=&\Psi(t),
\end{eqnarray}
where $Z(t)=A-BR^{-1}B'P(t)$, $\Phi(t)$ satisfies (\ref{ds4}), and $\Psi(t)=\Phi'(t)^{-1}$.
\end{corollary}

\emph{Proof:} The result follows directly from Lemma \ref{lem1}. \hfill $\blacksquare$

From Lemma \ref{lem1} and Corollary \ref{cor1}, we have obtained the numerical solution to (\ref{ds3}) and (\ref{ds4}). Accordingly,
the optimal solution (\ref{ds6}) can be numerically obtained, and then the detailed algorithm is given in Algorithm 1.
\begin{algorithm}[htb]
\caption{Centralized Algorithm}
\label{alg:Framwork}
    \textbf{Input}: Let the initial matrices be $P^{n}(T) = 0, n=1,2,3, \ldots$, $P^{0}(t)=0, t\in[0, T]$, and select a sufficiently small constant $\varepsilon> 0$.\\
     \textbf{Output}: The optimal states $x_{c}^{*}(t)$ and controller $u^{*}_c(t)$. \par
   \quad 1. For $n=1,2,3, \ldots, $ run the iteration (\ref{1}) and obtain $P^{n}(t)$. \par
   \quad 2. If $\|P^{n}(t)-P^{n-1}(t)\|<\varepsilon$, end the iteration, and derive $P^{n}(t)$, $Z^n(t)$ and  $\Phi^n(t)$. \par
    \quad3. Derive $\lambda_c^{*}$ by using  $$\lambda_c^{*}=\Big[\int_0^T\Phi^n(T,s)BR^{-1}B'\Psi^n(s,T)ds\Big]^{\dag}[\Phi^n(T,0)x_0-x_T].$$ \par
    \quad4. Solve the equation $$\dot{\beta}_c(t)=-Z^n(t)'\beta_c(t), \beta_c(T)=\lambda_c^{*},$$ and obtain $\beta_c(t)$. \par
    \quad5. Solve the equation $$\dot{x}_c^{*}(t)=Z^n(t)x^{*}(t)-BR^{-1}B'\beta_c(t), x_c^{*}(0)=x_0,$$ and obtain $x_c^{*}(t)$. \par
    \quad6. Let $u_c^{*}(t)=-R^{-1}B'P^n(t)x_c^{*}(t)-R^{-1}B'\beta_c(t)$ and output $u_c^{*}(t)$ and $x_{c}^{*}(t)$.
\end{algorithm}
\begin{theorem}\label{t2}
There exist an integer $n$ and a sufficient small constants $\epsilon_{1}> 0$, $\epsilon_{2}> 0$ such that output $x_c^{*}(t)$, $u_c^{*}(t)$ derived by Algorithm 1 satisfies
\begin{eqnarray}
\|x_c^{*}(t)-x^*(t)\|&\leq& \epsilon_{1},\label{c01}\\
\|u_c^{*}(t)-u^*(t)\|&\leq& \epsilon_{2},\label{c1}
\end{eqnarray}
where $x^{*}(t)$ and $u^*(t)$ are the optimal solution of equations (\ref{ds6}) and (\ref{ds9}) by assuming all the information of system parameters are available to make decision.
\end{theorem}

\emph{Proof:} From Lemma \ref{lem1}, there always exist an integer $n$ and a sufficient small constant $\varepsilon> 0$ such that $\|P^{n}(t)-P(t)\|\leq \varepsilon, t\in[0, T]$. This further implies that there exist sufficient small constants $\varepsilon_1> 0, \varepsilon_2> 0$ such that
\begin{eqnarray}
\|\lambda_c-\lambda^{*}\|\leq \varepsilon_1,\nonumber\\
 \|\beta_c(t)-\beta(t)\|\leq \varepsilon_2.\nonumber
\end{eqnarray}

Accordingly, there exists a sufficient small constant $\epsilon_{1}> 0$, $\epsilon_{2}> 0$ such that $x_c^{*}(t)$ and $u_c^{*}(t)$ obtained by running Algorithm 1 satisfies (\ref{c01}) and (\ref{c1}).
This completes the proof. \hfill $\blacksquare$

\begin{remark}
In Algorithm 1, solving the Riccati equation (\ref{1}) requires $Z^{n}(t)$ and $V^{n}(t)$, which depends on the global system matrices $A$ and $B$. However, this loses efficacy if each agent only has access to the local information as defined in Assumptions 1-4.
Something similar happens to $\lambda_c^{*}$, $x_c^{*}(t)$ and $\beta_c(t)$. This motivates the study of the distributed algorithm.
\end{remark}

\section{Main results}

In this section, we will present a distributed algorithm to derive optimal controller when only partial information is accessible to agent $i$ as provided in Assumptions 1-4 which is divided into three steps. The first step is to give the distributed iteration for $P(t)$ in (\ref{ds3}). The second step is to give the distributed iteration for optimal Lagrange multiplier $\lambda^*$ in (\ref{ds5}) and the distributed iteration for $\beta(t)$ in (\ref{ds8}). The third step is to give the distributed iteration for $x^*(t)$ in (\ref{ds9}).

\textit{\textbf{Step 1: }}Distributed iteration solution for $P(t)$

Firstly, we introduce the following distributed iteration of $P(t)$:
\begin{eqnarray}
\dot{P}_{i,k}^{n}(t)&=&-Z_{i,k}^{n}(t)'P_{i,k}^{n}(t)-P_{i,k}^{n}(t)
Z_{i,k}^{n}(t)-V_{i,k}^{n}(t),\label{31}
\end{eqnarray}
where $n=1,2,\cdots$, $k=1,2,\cdots$, $P_{i,k}^{0}(t)=0$, $P_{i,0}^{n}(t)=0$, $P_{i,k}^{n}(T)=0$,
\begin{align}
Z_{i,k}^{n}(t)=&Z_{i,k-1}^{n}(t)+\alpha_k[A_{i}-M_iM_i'P_{i,k-1}^{n-1}(t)-Z_{i,k-1}^{n}(t)]\nonumber\\
&+\frac{1}{\gamma}\sum_{j\in N_i}[Z_{j,k-1}^{n}(t)-Z_{i,k-1}^{n}(t)],\label{13}\\
V_{i,k}^{n}(t)=&V_{i,k-1}^{n}(t)+\alpha_k[Q_{i}+P_{i,k-1}^{n-1}(t)
M_iM_i'P_{i,k-1}^{n-1}(t)\nonumber\\
& -V_{i,k-1}^{n}(t)]+\frac{1}{\gamma}\sum_{j\in N_i}[V_{j,k-1}^{n}(t)\nonumber\\
&-V_{i,k-1}^{n}(t)],\label{14}
\end{align}
while $\alpha_k$ is a step size satisfying $\sum_{k=0}^{\infty}\alpha_k=\infty$, $\sum_{k=0}^{\infty}\alpha_k^2<\infty$, $\gamma>0$ is selected such that the matrix $I_N-\frac{1}{\gamma}\mathcal{L}-\frac{1}{N}1_N1_N'$ is Hurwitz, ensuring the convergence of the distributed iterations, and $Z_{i,0}^{n}(t)=0$, $V_{i,0}^{n}(t)=0$.

By solving (\ref{31}), the explicit solution of $P_{i,k}^{n}(t)$ is
\begin{eqnarray}
P_{i,k}^{n}(t)&=&\Psi_{i,k}^{n}(t)\big[\int_t^T\Phi_{i,k}^{n}(s)'
V_{i,k}^{n}(s)
\nonumber\\
&& \times\Phi_{i,k}^{n}(s)ds\big]\Psi_{i,k}^{n}(t)',\label{24}
\end{eqnarray}
where $\Phi_{i,k}^{n}(t)$ and $\Psi_{i,k}^{n}(t)$ satisfy the following iteration:
\begin{eqnarray}
\dot{\Phi}_{i,k}^{n}(t)&=&Z_{i,k}^{n}(t)\Phi_{i,k}^{n}(t), \label{18}\\
\dot{\Psi}_{i,k}^{n}(t)&=&-Z_{i,k}^{n}(t)'\Psi_{i,k}^{n}(t), \label{28}
\end{eqnarray}
with $\Phi_{i,k}^{n}(0)=I$ and $\Psi_{i,k}^{n}(0)=I$.

It is clear that the convergence of distributed iteration (\ref{31}) depends on the convergence of distributed iteration (\ref{13})-(\ref{14}) which is shown in the following result.

\begin{lemma}\label{lem2} It holds that
\begin{eqnarray}
\displaystyle\lim_{n\rightarrow \infty}\lim_{k\rightarrow \infty}\|Z_{i,k}^{n}(t)-Z(t)\|&=&0,\label{19}\\
\displaystyle\lim_{n\rightarrow \infty}\lim_{k\rightarrow \infty}\|V_{i,k}^{n}(t)-V(t)\|&=&0.\label{21}
\end{eqnarray}
\end{lemma}

\emph{Proof:} The detailed proof is given in Appendix A. \hfill $\blacksquare$

Then the following convergence result follows for the distributed iteration (\ref{31}), (\ref{18}) and (\ref{28}).

\begin{corollary}\label{cor2}
It holds that
\begin{eqnarray}
\displaystyle\lim_{n\rightarrow \infty}\lim_{k\rightarrow \infty}\|P_{i,k}^{n}(t)-P(t)\|&=&0,\\
\displaystyle\lim_{n\rightarrow \infty}\lim_{k\rightarrow \infty}\|\Phi_{i,k}^{n}(t)-\Phi(t)\|&=&0,\\
\displaystyle\lim_{n\rightarrow \infty}\lim_{k\rightarrow \infty}\|\Psi_{i,k}^{n}(t)-\Psi(t)\|&=&0.
\end{eqnarray}
\end{corollary}

\emph{Proof:} The result follows directly from Lemma \ref{lem2}. \hfill $\blacksquare$

For the convenience of marking, we make the denotations:
\begin{eqnarray}
&\displaystyle\lim_{n\rightarrow \infty}\lim_{k\rightarrow \infty}P_{i,k}^{n}(t)=P^{\infty}_{i}(t),\\
&\displaystyle\lim_{n\rightarrow \infty}\lim_{k\rightarrow \infty}\Phi_{i,k}^{n}(t)=\Phi^{\infty}_{i}(t), \\
&\displaystyle\lim_{n\rightarrow \infty}\lim_{k\rightarrow \infty}\Psi_{i,k}^{n}(t)=\Psi^{\infty}_{i}(t), \\
&\displaystyle\lim_{n\rightarrow \infty}\lim_{k\rightarrow \infty}Z_{i,k}^{n}(t)=Z^{\infty}_{i}(t), \\
&\displaystyle\lim_{n\rightarrow \infty}\lim_{k\rightarrow \infty}V_{i,k}^{n}(t)=V^{\infty}_{i}(t).
\end{eqnarray}

\textit{\textbf{Step 2: }}Distributed iteration solutions for $\lambda^*$ and $\beta(t)$

Secondly, to derive the distributed solution of $\lambda^{*}$, we denote
\begin{eqnarray}
\varrho= \int_0^T \Phi(s)W'W\Psi(s)ds,
\end{eqnarray}
where $W= R^{-\frac{1}{2}}B'$ and introduce the following distributed iteration for matrix $W$ as:
\begin{eqnarray}
W_{i,k}&=&W_{i,k-1}+\alpha_k[R^{-\frac{1}{2}}B_i'-W_{i,k-1}]\nonumber\\
&&+\frac{1}{\gamma}\sum_{j\in N_i}[W_{j,k-1}-W_{i,k-1}],\label{29}
\end{eqnarray}
with $W_{i,0}=0$, whose convergence is given below.

\begin{lemma}\label{lem3}It holds that
\begin{eqnarray}
W_{i}^{\infty}\triangleq \displaystyle\lim_{k\rightarrow \infty}W_{i,k}=W.
\end{eqnarray}
\end{lemma}

\emph{Proof:} The proof is similar to Lemma 2 and is omitted here.  \hfill $\blacksquare$

Recalling (\ref{ds5}) where $\lambda^{*}$ depends not only on the global $W$ but also on the global initial state $x(0)$ and terminal state $x(T)$, we then propose the distributed iteration of $\lambda^{*}$ by
\begin{eqnarray}
\lambda_{i,k}^{*}&=&\lambda_{i,k-1}^{*}+\alpha_k\{(\varrho_{i}^{\infty})^{\dag}[\Phi^{\infty}_{i}(T,0)x_{i,0}-x_{i,T}]
\nonumber\\
&&-\lambda_{i,k-1}^{*}\}+\frac{1}{\gamma}\sum_{j\in N_i}[\lambda_{j,k-1}^{*}-\lambda_{i,k-1}^{*}],\label{30}
\end{eqnarray}
with $\lambda_{i,0}^{*}=0$ and
\begin{eqnarray}
\varrho_{i}^{\infty}= \int_0^T \Phi^{\infty}_{i}(s){W_{i}^{\infty}}'W_{i}^{\infty}\Psi^{\infty}_{i}(s)ds.
\end{eqnarray}
The convergence of this iteration in the following result.

\begin{lemma}\label{lem4}
It holds that
\begin{eqnarray}
\lambda_{i}^{\infty,*}\triangleq\displaystyle\lim_{k\rightarrow \infty}\lambda_{i,k}^{*}=\lambda^{*}.
\end{eqnarray}
\end{lemma}

\emph{Proof:} The proof is similar to Lemma 2 and is omitted here. \hfill $\blacksquare$

Accordingly, the distributed solution for $\beta(t)$ is defined as follows:
\begin{eqnarray}
\dot{\beta}_{i}^{\infty}(t)&=&-Z^{\infty}_{i}(t)'\beta_{i}^{\infty}(t),\label{36}\\
\beta_{i}^{\infty}(T)&=&\lambda_{i}^{\infty,*}.\nonumber
\end{eqnarray}

\textit{\textbf{Step 3:}} Distributed iteration solution for $x^{*}(t)$

Thirdly, we give the distributed solution of $x^{*}(t)$. In fact, by solving (\ref{ds5}), the explicit solution of $x^{*}(t)$ is given by
\begin{eqnarray}
x^{*}(t)=\Phi(t,0)x_{0}-\int_0^t\Phi(t,s)W'W\beta(s) ds,\label{35}
\end{eqnarray}
which relies on the global initial state $x_{0}$. To avoid the using of global information, we give the distributed iteration as follows:
\begin{eqnarray}
x_{i,k}^{*}(t)&=&x_{i,k-1}^{*}(t)+\alpha_k[\Phi_{i}^{\infty}(t,0)x_{i,0}-\int_0^t
\Phi_{i}^{\infty}(t,s)\nonumber\\
&&\times {W_{i}^{\infty}}'W_{i}^{\infty}\beta_{i}^{\infty}(s) ds-x_{i,k-1}^{*}(t)]\nonumber\\
&&
+\frac{1}{\gamma}\sum_{j\in N_i}[x_{j,k-1}^{*}(t)-x_{i,k-1}^{*}(t)], \label{15}
\end{eqnarray}
with $x_{i,0}^{*}(t)=0$ and derive the following convergence result.

\begin{lemma}\label{lem5} It holds that
\begin{eqnarray}
\displaystyle\lim_{k\rightarrow \infty}x_{i,k}^{*}(t)=x^{*}(t).
\end{eqnarray}
\end{lemma}

\emph{Proof:} The proof is similar to Lemma 2 and is omitted here.  \hfill $\blacksquare$

From Lemmas \ref{lem2}-\ref{lem5} and Corollary \ref{cor2}, we have established the distributed solutions to  (\ref{ds3}) and (\ref{ds4}). This enables the distributed design of the optimal controller (\ref{ds6}) as shown in Algorithm 2.
\begin{algorithm}[htb]
\caption{Distributed Algorithm}
\label{alg:Distributed}
\textbf{Input}: Let the initial matrices be $Z_{i,0}^{n}(t)=0$, $V_{i,0}^{n}(t)=0$, $P_{i,k}^{0}(t)=0$, $W_{i,0}=0$, $\lambda_{i,0}^{*}=0$, and $x_{i,0}^{*}(t)=0$ and select a sufficiently small constant $\varsigma> 0$. \par
\textbf{Output}: The optimal states $x_{i,w}^{*}(t)$ and distributed controller $u_i^{*}(t)$, $i=1, 2, \ldots, N$. \par
\begin{enumerate}
\item[1.] For $n=1,2,\cdots$,
\begin{itemize}
\item For $k=1,2,3, \cdots$, run the iteration (\ref{13})-(\ref{14}) and obtain $Z_{i,k}^{n}(t)$, $V_{i,k}^{n}(t)$. If $\|Z_{i,k}^{n}(t)-Z_{i,k-1}^{n}(t)\|<\varsigma$ and $\|V_{i,k}^{n}(t)-V_{i,k-1}^{n}(t)\|<\varsigma$, end the iteration, and derive $Z_{i,k}^{n}(t)$ and $V_{i,k}^{n}(t)$. \par
  \item Solve the equation (\ref{31}) and obtain $P_{i,k}^{n}(t)$. If $\|P_{i,k}^{n}(t)-P_{i,k}^{n-1}(t)\|<\varsigma$, end the iteration, and derive $P^{n}_{i,k}(t)$, $\Phi^{n}_{i,k}(t)$, $\Psi^{n}_{i,k}(t)$. \par
\end{itemize}
\item[2.] For $\varpi=1,2,\cdots$, run the iteration
\begin{eqnarray*}
W_{i,\varpi}&=&W_{i,\varpi-1}+\alpha_\varpi[R^{-\frac{1}{2}}B_i'-W_{i,\varpi-1}]\\
&&+\frac{1}{\gamma}\sum_{j\in N_i}[W_{j,\varpi-1}-W_{i,\varpi-1}].
\end{eqnarray*}
\begin{itemize}
\item If $\|W_{i,\varpi}-W_{i,\varpi-1}\|<\varsigma$, end the iteration, and derive $W_{i,\varpi}$. \par
\end{itemize}
\item[3.] For $q=1,2,\cdots$, run the iteration
\begin{eqnarray*}
\lambda_{i,q}^{*}&=&\lambda_{i,q-1}^{*}+\alpha_q\bigg\{\Big[ \int_0^T \Phi_{i,k}^{n}(s)W_{i,\varpi}'W_{i,\varpi}
\nonumber\\
&&\times\Psi_{i,k}^{n}(s)ds\Big]^{\dag}[\Phi_{i,k}^{n}(T,0)x_{i,0}-x_{i,T}]
\nonumber\\
&&-\lambda_{i,q-1}^{*}\}+\frac{1}{\gamma}\sum_{j\in N_i}[\lambda_{j,q-1}^{*}-\lambda_{i,q-1}^{*}].
\end{eqnarray*}
\begin{itemize}
\item If $\|\lambda_{i,q}^{*}-\lambda_{i,q-1}^{*}\|<\varsigma$, end the iteration, and derive $\lambda_{i,q}^{*}$. \par
\end{itemize}
\item[4.] Solve the equation
\begin{eqnarray*}
\dot{\beta}_i(t)= -Z_{i,k}^{n}(t)'\beta_i(t), ~\beta_i(T)=\lambda_{i,q}^{*},
\end{eqnarray*}
and obtain $\beta_i(t)$. \par
\item[5.] For $w=1,2,\cdots$, run the iteration
\begin{eqnarray*}
x_{i,w}^{*}(t)&=&x_{i,w-1}^{*}(t)+\alpha_w[\Phi_{i,k}^{n}(t,0)x_{i,0}-\int_0^t
\Phi_{i,k}^{n}(t,s)\nonumber\\
&&\times W_{i,\varpi}'W_{i,\varpi}\beta_{i}(s) ds-x_{i,w-1}^{*}(t)]\nonumber\\
&&
+\frac{1}{\gamma}\sum_{j\in N_i}[x_{j,w-1}^{*}(t)-x_{i,w-1}^{*}(t)].
\end{eqnarray*}
\begin{itemize}
\item If $\|x_{i,w}^{*}(t)-x_{i,w-1}^{*}(t)\|<\varsigma$, end the iteration, and derive $x_{i,w}^{*}(t)$. \par
\end{itemize}
\item[6.] Let $u_i^{*}(t)=-R^{-1/2}W_{i,\varpi}\big[P_{i,k}^n(t)x_{i,w}^{*}(t)
    +\beta_i(t)\big]$ and output $u_i^{*}(t)$ and $x_{i,w}^{*}(t)$.
\end{enumerate}
\end{algorithm}

\begin{theorem}\label{t3}
There exist integers $n,k,\varpi,q,w$ and sufficiently small constants $\psi_{1}>0$, $\psi_{2}>0$ such that output $x_{i,w}^{*}$ and $u_i^{*}(t)$ derived by Algorithm 2 satisfy:
\begin{eqnarray}
\|x_{i,w}^*(t)-x^*(t)\|\leq \psi_{1},\label{cc2}\\
\|u_i^*(t)-u^*(t)\|\leq \psi_{2},\label{c2}
\end{eqnarray}
where $x^{*}(t)$ and $u^*(t)$ are the optimal solution of equations (\ref{ds6}) and (\ref{ds9}) by assuming all the information of system parameters are available to make decision.
\end{theorem}

\emph{Proof.} From Lemmas \ref{lem1}-\ref{lem2} and Corollary \ref{cor2}, there always exists integer $n$, $k$ and a sufficient small constant $\sigma_1> 0$ such that $\|P_{i,k}^{n}(t)-P(t)\|\leq \sigma_1, t\in[0, T]$. Together with Lemmas \ref{lem3}-\ref{lem4}, this further implies that there exist iteration numbers $\varpi,q$ and sufficient small constants $\sigma_{i}> 0, i=2,\cdots,4$ such that
\begin{eqnarray}
\|W_{i,\varpi}-W\|\leq \sigma_2,\nonumber\\
\|\lambda_{i,q}^{*}-\lambda^{*}\|\leq \sigma_3,\nonumber\\
 \|\beta_i(t)-\beta(t)\|\leq \sigma_4.\nonumber
\end{eqnarray}

Accordingly, by using Lemma \ref{lem5}, there exist iteration numbers $w$ and sufficient small constants $\psi_1> 0$, $\psi_{2}> 0$ such that the outputs $x_{i,w}^*(t)$ and $u_i^{*}(t)$ obtained by running Algorithm 2 satisfies (\ref{cc2}) and (\ref{c2}). This completes the proof.\hfill $\blacksquare$

\section{Application to the Optimal Consensus Control Problem}

In this section, we will utilize Algorithm 2 to solve the consensus problem in multi-agent systems, which has broad applications in fields such as unmanned aerial vehicle formation, sensor networks, smart grid dispatch, and distributed robotic collaboration. Specifically, consider the following heterogeneous multi-agent system consisting of $N$ agents with dynamics described by:
\begin{align}
\dot{x}_{i}(t)=A_{i}x_{i}(t)+B_{i}u_{i}(t), \label{cc1}
\end{align}
where $x_{i}(t)\in R^{n}$, $u_{i}(t)\in R^{m}$, $i=1,..,N$. The initial state is given by $x_{i}(0)=x_{i,0}$.

To apply Algorithm 2 and achieve consensus, we define the following performance index:
\begin{align}
J&=\int_{0}^{\infty}\Big[\sum_{i=1}^{N}\sum_{j\in N_i}(x_i(t)-x_j(t))'Q_{ij}(x_i(t)-x_j(t)) \nonumber \\
&\quad\quad\quad~~+u'(t)Ru(t)\Big]dt\nonumber \\
&=\int_{0}^{\infty}[x'(t)\tilde{Q}x(t)+u'(t)Ru(t)]dt,\label{cc2}
\end{align}
where $x(t)=[x_{1}'(t)~\cdots~x_{N}'(t)]'$, $u(t)=[u_{1}'(t)~\cdots~u_{N}'(t)]'$,
$R=diag\{R_{1}, R_{2},\ldots, R_{N}\}$ is a positive-definite matrix with compatible dimension, $Q_{ij}$, $j \in N_{i}$ is a positive definite matrix with compatible dimension, and
\begin{eqnarray}
\tilde{Q}&=\begin{bmatrix}
\tilde{Q}_{11}& \tilde{Q}_{12}& \cdots &\tilde{Q}_{1N} \nonumber\\
\vdots& \vdots& \cdots& \vdots \nonumber\\
\tilde{Q}_{N1}&\tilde{Q}_{N2}& \cdots &\tilde{Q}_{NN}\nonumber
\end{bmatrix},
\end{eqnarray}
with
\begin{align}
\tilde{Q}_{ii}&=\sum_{j=1}^{N}c_{ij}(Q_{ij}+Q_{ji}), \nonumber\\
\tilde{Q}_{ij}&=-c_{ij}(Q_{ij}+Q_{ji}), \nonumber
\end{align}
while $c_{ij} = 1$ if $j \in N_{i}$, otherwise $c_{ij} = 0$. In particular, each agent $i$ has access to the $i$-th block row of $\tilde{Q}$, i.e., $\tilde{Q}_i= L_i\tilde{Q}$ with $L_i= diag\{0,\cdots,0, NI,0,\cdots, 0\}$.

The objective of this section is to design a distributed controller by using its own and neighboring information of the system (\ref{cc1}) to minimize the cost function (\ref{cc2}).

\begin{remark}
In the literature, a commonly used protocol for multi-agent systems \cite{37,38,39} is:
\begin{align}
u_i(t)=K\sum_{j\in N_i}(x_{j}(t)-x_{i}(t)),\label{cc15}
\end{align}
which ensures that the system (\ref{cc1}) achieves state consensus, that is, for any initial state $x_i(0)$,
 \begin{align}
\lim_{t\rightarrow \infty}\|x_i(t)-x_j(t)\|=0,~~i,j=1,2,\cdots,N.\label{cc14}
\end{align}
However, such protocols do not consider the optimality, and the optimality of some specific performance index, which is essential in practical scenarios such as smart grids \cite{16} and unmanned ground vehicle (UGV) swarm coordination \cite{32}. In this section, we present an optimal consensus protocol using Algorithm 2.
\end{remark}

To address this, we first present the centralized optimal controller, which requires information from all agents.
\begin{lemma}\label{lem6}
Assume that $(A_i, B_i)$ is stabilizable and $(A, \tilde{Q}^{1/2})$ is observable, where $A=diag\{A_1,\cdots,A_N\}$. Then, the centralized optimal controller is given by
\begin{align}
u^{\diamond}(t) &= -R^{-1} B' P x^{\diamond}(t), \label{cc3}
\end{align}
where $B=diag\{B_1,\cdots, B_N\}$, $P$ is the unique positive-definite solution to the algebraic Riccati equation:
\begin{align}
A' P + P A + \tilde{Q} - P B R^{-1} B' P &= 0, \label{cc4}
\end{align}
and $x^{\diamond}(t)$ satisfies:
\begin{align}
\dot{x}^{\diamond}(t) &= [A - B R^{-1} B' P] x^{\diamond}(t), \quad x^{\diamond}(0) = x(0). \label{cc5}
\end{align}
\end{lemma}

\emph{Proof.} Since $(A_i, B_i)$ is stabilizable, the block diagonal $(A, B)$ is stabilizable. Combined with the observability of $(A, \tilde{Q}^{1/2})$, the algebraic Riccati equation \eqref{cc4} admits a unique positive-definite solution $P$. The optimal controller \eqref{cc3} and the corresponding closed-loop system \eqref{cc5} then follow directly from standard LQR theory \cite{34}. This completes the proof. \hfill $\blacksquare$

Next, we derive the distributed solution by using Algorithm 2. To this end, we define
\begin{align}
&\tilde{A}_{i} = diag\{0,\cdots,0,NA_{i},0,\cdots,0\},\nonumber\\
&\tilde{B}_{i} = diag\{0,\cdots,0,NB_{i},0,\cdots,0\},\nonumber\\
&\tilde{R}_i = diag\{0,\cdots, (1/N)R_{i}^{-1},\cdots, 0\},\nonumber
\end{align}
and introduce the distributed iteration equation for $P$:
\begin{align}
0&=(\bar{Z}_{i,k}^{n})'P_{i,k}^{n}+P_{i,k}^{n}\bar{Z}_{i,k}^{n}+\bar{V}_{i,k}^{n},\label{cc6}\\
\bar{Z}_{i,k}^{n}&=\bar{Z}_{i,k-1}^{n}+\alpha_k[\tilde{A}_{i}-\tilde{B}_{i}\tilde{R}_{i}\tilde{B}_{i}'P_{i,k-1}^{n-1}-\bar{Z}_{i,k-1}^{n}]\nonumber\\
&\quad+\frac{1}{\gamma}\sum_{j\in N_i}[\bar{Z}_{j,k-1}^{n}-\bar{Z}_{i,k-1}^{n}],\label{cc7}\\
\bar{V}_{i,k}^{n}&=\bar{V}_{i,k-1}^{n}+\alpha_k[\tilde{Q}_{i}+P_{i,k-1}^{n-1}
\tilde{B}_{i}\tilde{R}_{i}\tilde{B}_{i}'P_{i,k-1}^{n-1}\nonumber\\
&\quad -\bar{V}_{i,k-1}^{n}]+\frac{1}{\gamma}\sum_{j\in N_i}[\bar{V}_{j,k-1}^{n}-\bar{V}_{i,k-1}^{n}],\label{cc8}
\end{align}
with $\bar{Z}_{i,0}^{n}=0$, $P_{i,k}^{0}=0$, $P_{i,0}^{n}=0$, $\bar{V}_{i,0}^{n}=0$, and the distributed iteration equation for $x^{\diamond}(t)$:
\begin{align}
x_{i,k}^{\diamond}(t) &= x_{i,k-1}^{\diamond}(t) + \alpha_k \big[e^{\bar{Z}_{i,\infty}t} \bar{x}_{i}(0) -x_{i,k-1}^{\diamond}(t)\big] \nonumber \\
& \quad + \frac{1}{\gamma}\sum_{j \in \mathcal{N}_i} [x_{j,k-1}^{\diamond}(t) - x_{i,k-1}^{\diamond}(t)], \label{cc9}
\end{align}
with $\bar{x}_{i}(0)= [0~\cdots~0~Nx_{i,0}'~0~\cdots~ 0]'$,  $x_{i,0}^{\diamond}(t)=0$ and $\bar{Z}_{i,\infty}$ is the convergent value of (\ref{cc7}).

Then, we present the distributed design of the controller (\ref{cc3}) as follows.
\begin{theorem}\label{t4}
Under the stabilizability of $(A_i,B_i)$ and observability of $(A,\tilde{Q}^{1/2})$, let $\alpha_k$ satisfy $\sum_{k=0}^\infty \alpha_k = \infty$, $\sum_{k=0}^\infty \alpha_k^2 < \infty$, and choose $\gamma > 0$ such that $I_N - \frac{1}{\gamma}\mathcal{L} - \frac{1}{N}\mathbf{1}_N\mathbf{1}_N'$ is Hurwitz, then it holds that:
\begin{align*}
P_{i,\infty} &\triangleq \lim_{n\to\infty} \lim_{k\to\infty} P_{i,k}^{n}=P,~~ \\
\bar{Z}_{i,\infty}& \triangleq \lim_{n\to\infty} \lim_{k\to\infty} \bar{Z}_{i,k}^{n}=A-BR^{-1}B'P,\\
x_{i,\infty}^{\diamond}(t) &\triangleq \lim_{k\to\infty} x_{i,k}^{\diamond}(t)=x^{\diamond}(t).
\end{align*}

In this case, the distributed optimal controller for each agent $i$ is given by
\begin{align}
u_i(t) = -[0~ \ldots~  0~  R_i^{-1}B_{i}'~  0~  \ldots~  0] P_{i,\infty} x_{i,\infty}^{\diamond}(t), \label{cc10}
\end{align}
which satisfies
\begin{align}
\|u(t) - u^{\diamond}(t)\| \leq \varepsilon, \label{cc11}
\end{align}
and ensures that the closed-loop system (\ref{cc1}) achieves consensus.
\end{theorem}

\emph{Proof.} The proof is divided into two parts. First, we prove that (\ref{cc11}) holds. Similar to the proof of Lemma \ref{lem2}, there exist integers $\bar{n}\geq n,\bar{k}\geq k$ and sufficiently small constants $\varepsilon_{1}>0$, $\varepsilon_{2}>0$ such that
\begin{align}
\|P_{i,\bar{k}}^{\bar{n}}-P\|&\leq \varepsilon_{1},\label{cc12}\\
\|x_{i,\bar{k}}^{\diamond}(t)-x^{\diamond}(t)\|&\leq \varepsilon_{2}.\label{cc13}
\end{align}
This implies that the distributed solutions $P_{i,k}^{n}$ and $x_{i,k}^{\diamond}(t)$ can approximate their global solution $P$ and $x^{\diamond}(t)$. Then, combining the fact that $u(t)=[u_{1}'(t)~\cdots~u_{N}'(t)]'$ and
\begin{eqnarray}
 \begin{bmatrix}
R_1^{-1}B_{1}'& 0& \cdots &0 \nonumber\\
0& R_2^{-1}B_{2}'& \cdots &0 \nonumber\\
\vdots&  \vdots & \ddots& \vdots \nonumber\\
0&0& \cdots &R_N^{-1}B_{N}'\nonumber
\end{bmatrix}=R^{-1}B',
\end{eqnarray}
we can design the distributed controller (\ref{cc10}) from (\ref{cc3}) and obtain (\ref{cc11}) holds.

Second, we show that the system achieves consensus. Since the problem has a solution and $J<\infty$, then it follows from the properties of infinite integrals that
 \begin{align}
\lim_{t\rightarrow \infty}(x_i(t)-x_j(t))'Q_{ij}(x_i(t)-x_j(t))=0.\label{cc16}
\end{align}
Combined with $Q_{ij}>0$, the necessary condition for the term (\ref{cc16}) to converge to zero over the infinite horizon is
\begin{align*}
\lim_{t\rightarrow \infty}\|x_i(t)-x_j(t)\|=0,~~i,j=1,2,\cdots,N.\label{cc14}
\end{align*} This completes the proof.\hfill $\blacksquare$
\section{Simulation Example}
In this section, two examples are provided to show the effectiveness of the proposed Algorithms.

Example 1. Consider system (\ref{ds1}) and the cost function (\ref{ds2}) with parameters given by
\begin{eqnarray}
&&A =
\left[
\begin{array}{cc}
1~&0\\
0~&1
\end{array}\right], B =
\left[
\begin{array}{cc}
1\\
1
\end{array}\right], Q=I, R=1, \\
&&x_{0} =\left[\begin{array}{cc}
4\\
4
\end{array}\right],
x_{T} =\left[\begin{array}{cc}
0\\
0
\end{array}\right].
\end{eqnarray}

Moreover, 4 agents are involved in the sensor network and the topology graph is shown in Fig. 1.
\begin{figure}
\begin{center}
  \includegraphics[width=5.2cm,height=4.6cm]{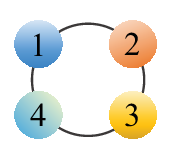}\\
  {\footnotesize Fig. 1. Communication topology among agents.}
  \end{center}
\end{figure}

Next, we consider the distributed case where each agent only has access to partial system information. Based on Assumptions 1-4, the parameters for each agent are given by
\begin{eqnarray*}
&&A_1 =\left[\begin{array}{cc}1.5~&0\\0~&0.5\end{array}\right],
A_2 =\left[\begin{array}{cc}0.5~&0\\0~&1.2\end{array}\right],\\
&& A_3 =\left[\begin{array}{cc}1.3~&0\\0~&1.8\end{array}\right],
A_4 =\left[\begin{array}{cc}0.7~&0\\0~&0.5\end{array}\right],\\
&&
B_1 =\left[\begin{array}{cc}1.5\\0.5\end{array}\right],
B_2 =\left[\begin{array}{cc}0.5\\1.2\end{array}\right],
B_3 =\left[\begin{array}{cc}0.3\\1.8\end{array}\right],
\\
&&B_4 =\left[\begin{array}{cc}1.7\\0.5\end{array}\right],
Q_1 =\left[\begin{array}{cc}1.2~&0\\0~&1\end{array}\right],
Q_2 =\left[\begin{array}{cc}1~&0\\0~&0.8\end{array}\right],\\
&&
Q_3 =\left[\begin{array}{cc}1~&0\\0~&1.2\end{array}\right],
Q_4 =\left[\begin{array}{cc}0.8~&0\\0~&1\end{array}\right],
R=1,\\
&&
x_{1,0} =\left[\begin{array}{cc}3.5\\4\end{array}\right],
x_{2,0} =\left[\begin{array}{cc}4.5\\3.5\end{array}\right],
x_{3,0} =\left[\begin{array}{cc}2.5\\4.5\end{array}\right],\\
&&
x_{4,0} =\left[\begin{array}{cc}5.5\\4\end{array}\right],
x_{i,T} =\left[\begin{array}{cc}0\\0\end{array}\right],
M_i =\left[\begin{array}{cc}1~&0\\0~&1\end{array}\right],\\
&&i=1,\cdots,4.
\end{eqnarray*}

By running Algorithm 2 with $\alpha_{k}=\frac{1}{k}$, $\gamma=2.5$, $\varsigma=10^{-3}$, $n=20$, $k=\varpi=q=w=200$, we obtain $P_{i,k}^{n}(t)$, $W_{i,\varpi}$, $\lambda_{i,q}^{*}$, $\beta_{i}(t)$, $x_{i,w}^{*}(t)$ and $u_{i}^{*}(t)$, and the corresponding trajectories are shown in Figs. 2-7. It is obvious that the distributed solutions converge to the centralized solutions from the Figs. 2-7.
\begin{figure}
\begin{center}
  \includegraphics[width=8.2cm,height=6.3cm]{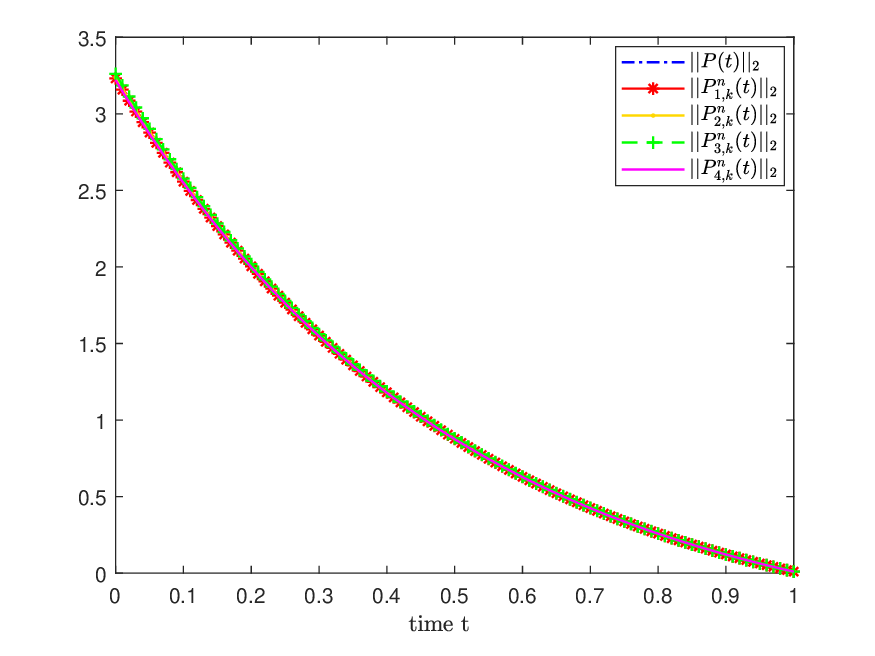}\\
  {\footnotesize Fig. 2. The trajectories of $\|P_{i,k}^{n}(t)\|_2$.}
  \end{center}
\end{figure}
\begin{figure}
\begin{center}
  \includegraphics[width=8.2cm,height=6.3cm]{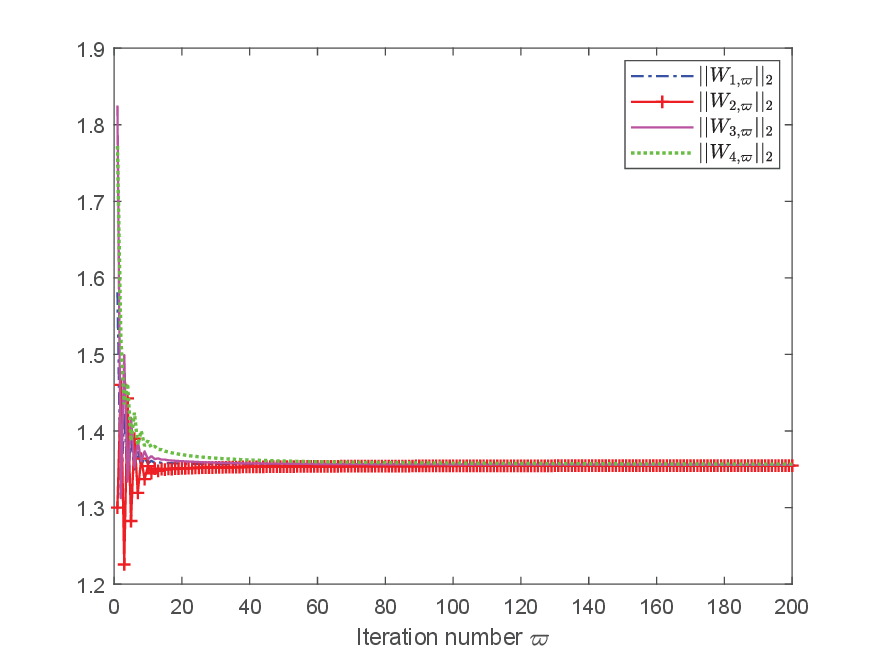}\\
  {\footnotesize Fig. 3. The trajectories of $\|W_{i,\varpi}\|_2$.}
  \end{center}
\end{figure}
\begin{figure}
\begin{center}
  \includegraphics[width=8.2cm,height=6.3cm]{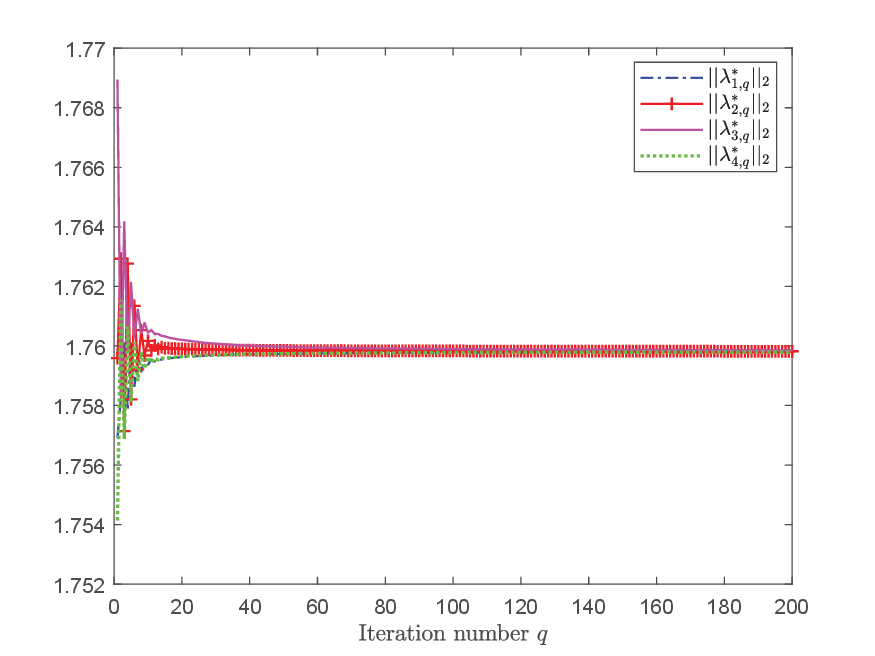}\\
  {\footnotesize Fig. 4. The trajectories of $\|\lambda_{i,q}^{*}\|_2$.}
  \end{center}
\end{figure}
\begin{figure}
\begin{center}
  \includegraphics[width=8.2cm,height=6.3cm]{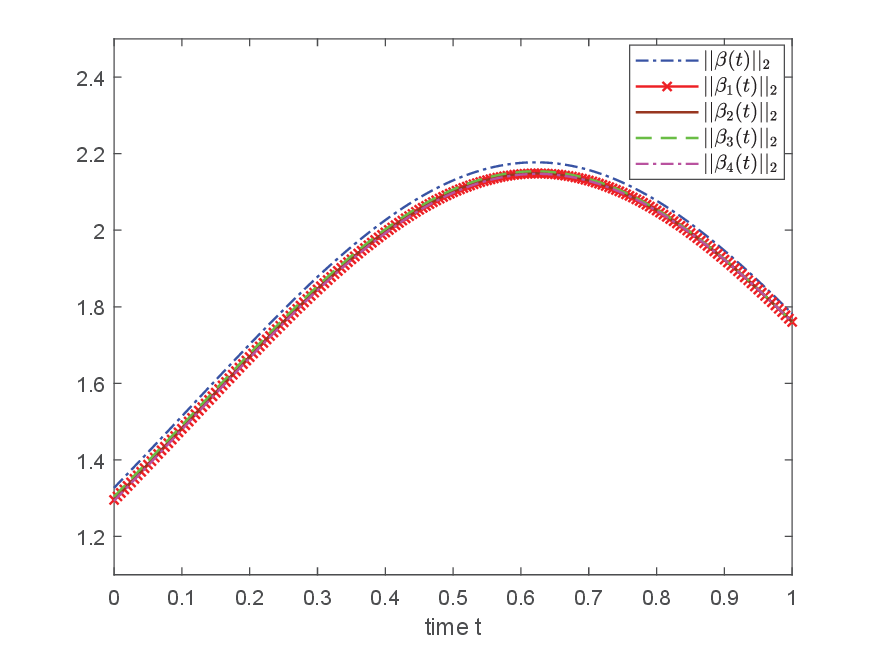}\\
  {\footnotesize Fig. 5. The trajectories of $\|\beta_{i}(t)\|_2$.}
  \end{center}
\end{figure}
\begin{figure}
\begin{center}
  \includegraphics[width=8.2cm,height=6.3cm]{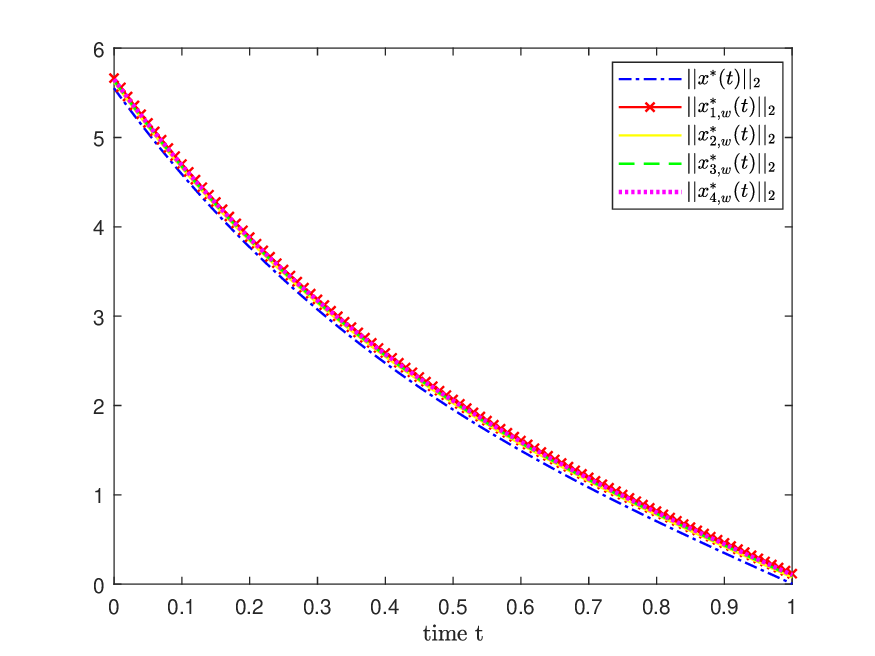}\\
  {\footnotesize Fig. 6. The trajectories of $\|x_{i,w}^{*}(t)\|_2$.}
  \end{center}
\end{figure}
\begin{figure}
\begin{center}
  \includegraphics[width=8.2cm,height=6.3cm]{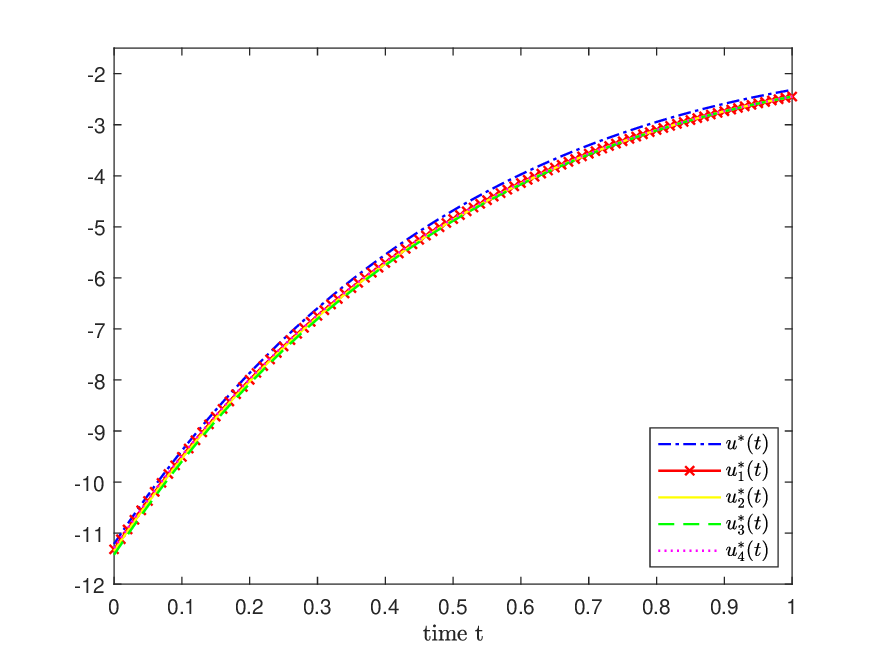}\\
  {\footnotesize Fig. 7. The trajectories of $u_{i}^{*}(t)$.}
  \end{center}
\end{figure}

Example 2. In the example, we will apply Theorem 4 to solve the consensus problem of a heterogeneous UGV system. Specifically, the objective is to enable five UGVs to achieve position consensus, where the dynamics of each UGV are given by \cite{35}:
\begin{equation}
\left[\begin{array}{c}
\dot{q}_{i}(t) \\
\dot{v}_{i}(t)
\end{array}\right] =
\left[\begin{array}{cc}
0 & I \\
0 & -\frac{C_i}{D_i}I
\end{array}\right]
\left[\begin{array}{c}
q_{i}(t) \\
v_{i}(t)
\end{array}\right] +
\left[\begin{array}{c}
0 \\
\frac{1}{D_i}I
\end{array}\right] F_i(t),
\label{si1}
\end{equation}
while $q_{i}(t)\in R^{2}$, $v_{i}(t)\in R^{2}$, $i = 1,\cdots,5$, denote the position and velocity vectors of the $i$th UGV, $F_{i}(t)\in R^{2}$ denotes the force applied to the UGV, $D_{i}$ denotes the mechanical mass and $C_{i}$ denotes the translational friction coefficient. The system parameters are listed in Table I, and the communication topology among UGVs is illustrated in Fig. 8.
\begin{figure}
\begin{center}
  \includegraphics[width=5.3cm,height=4.8cm]{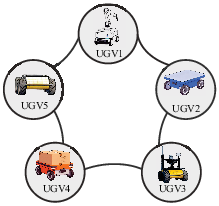}\\
  {\footnotesize Fig. 8. Communication topology among heterogeneous UGVs.}
  \end{center}
\end{figure}
\begin{table}[t]
\centering
\caption{Parameters of the UGV system.}
\label{table1}
\begin{tabular}{l|l|l|l|l|l}

\hline  
Parameters & UGV1 & UGV2  & UGV3 & UGV4 & UGV5 \\
\hline 
$C_i$ & $0.6$ & $0.8$ & $1.2$ & $1$& $0.4$ \\
\hline 
$D_i$ & $5$ & $4$& $6$ & $4$& $3$  \\
\hline 
\end{tabular}
\end{table}

By defining the state variables as $x_{i}(t)=[q_{i}(t)'~~v_{i}(t)']'$ and the control input as $u_{i}(t)=F_{i}(t)$, then the UGV system (\ref{si1}) can be formulated as (\ref{cc1}) with
\begin{eqnarray*}
&&A_{i} =
\left[\begin{array}{cc}
0 & I \\
0 & -\frac{C_i}{D_i}I
\end{array}\right], ~~
B_{i}=
\left[\begin{array}{c}
0 \\
\frac{1}{D_i}I
\end{array}\right].
\end{eqnarray*}

Let the parameters in cost function (\ref{cc2}) be given by
\begin{eqnarray*}
&&Q_{11} = Q_{22}=Q_{33}=Q_{44}=Q_{55}=0,\\
&&Q_{12}=Q_{15}=Q_{23} = Q_{34} = Q_{45}=I,\\
&&Q_{21}=Q_{51}=Q_{32} = Q_{43} = Q_{54}=I,\\
&&Q_{13}=Q_{14}=Q_{24} = Q_{25} = Q_{35}=0,\\
&&Q_{31}=Q_{41}=Q_{42} = Q_{52} = Q_{53}=0,\\
&&R_1=R_2=R_3=R_4=R_5=I.
\end{eqnarray*}

In particular, the initial positions and velocities of the UGVs are given by
\begin{eqnarray*}
&&x_{1}(0) =\left[\begin{array}{cccc}2\\6\\1\\1\end{array}\right],
x_{2}(0) =\left[\begin{array}{cccc}3\\3\\2\\2\end{array}\right],
x_{3}(0) =\left[\begin{array}{cccc}2\\2\\1\\2\end{array}\right],\\
&&x_{4}(0)=\left[\begin{array}{cccc}1\\2\\1.5\\1.5\end{array}\right],
x_{5}(0) =\left[\begin{array}{cccc}1\\4\\2\\1\end{array}\right].
\end{eqnarray*}

By selecting step size $\alpha_k = 1/k$, $\gamma=2.5$, $\varepsilon=\varepsilon_1=\varepsilon_2 = 10^{-4}$, we execute the distributed iterative equations (\ref{cc6})-(\ref{cc9}). Then it follows that:
\begin{itemize}
\item The distributed solutions $\|P_{i,\bar{k}}^{\bar{n}}\|_2$ with $\bar{n}=20$, $\bar{k}=900$ are:
  \begin{align*}
  &\|P_{1,\bar{k}}^{\bar{n}}\|_2=33.5490, ~\|P_{2,\bar{k}}^{\bar{n}}\|_2=33.5134, \\ &\|P_{3,\bar{k}}^{\bar{n}}\|_2=33.5116, ~\|P_{4,\bar{k}}^{\bar{\bar{n}}}\|_2=33.5111, \\ &\|P_{5,\bar{k}}^{n}\|_2=33.5475,
  \end{align*}
  which approximate the centralized solution $\|P\|_2=33.5591$ solved by (\ref{cc4}).
  \item The trajectories of $x_{i,\bar{k}}^{\diamond}(t)$ for all UGVs are shown in Fig. 9. It is obvious that $x_{i,k}^{\diamond}(t)$ converges to the centralized solution $x^{\diamond}(t)$ solved by (\ref{cc5}).
  \item The trajectory of $u(t)$, composed of all distributed controllers $u_i(t)$, is shown in Fig. 10, which demonstrates that $u(t)$ converges to the centralized optimal controller $u^{\diamond}(t)$ designed in (\ref{cc3}).
  \item The position and velocity trajectories of the five UGVs are shown in Figs. 11 and 12, respectively, which demonstrate that all UGVs achieve position consensus, and their velocities eventually converge to zero.
\end{itemize}

In particular, we conduct a comparative simulation with the widely used classical consensus protocol \cite{38} under the same initial conditions and system parameters $C_i$, $D_i$. The results are shown in Table II, from which it is obvious that the performance index $J$ of the proposed method is smaller than that obtained by using the controllers in \cite{38}.

\begin{table}[h]
\centering
\caption{Performance index comparison between the proposed method and [38] under different $(C_i, D_i)$.}
\label{table2}
\begin{tabular}{l|l|l|l}
\hline  
Case & $C_i, D_i$ & $J$ (Proposed Method) & $J$ ([38])   \\
\hline 
Case 1 & (0.6, 5)& 82.7573 &  125.1136   \\
\hline 
Case 2 & (0.8, 4)& 76.2762 &   107.3754   \\
\hline 
Case 3 & (1.2, 6)& 89.3999 &   135.6383   \\
\hline 
Case 4 & (1, 4)& 76.4114 &   106.4391   \\
\hline 
Case 5 & (0.4, 5)& 82.7004 &   127.5681  \\
\hline 
\end{tabular}
\end{table}
\begin{figure}
\begin{center}
  \includegraphics[width=8cm,height=6.3cm]{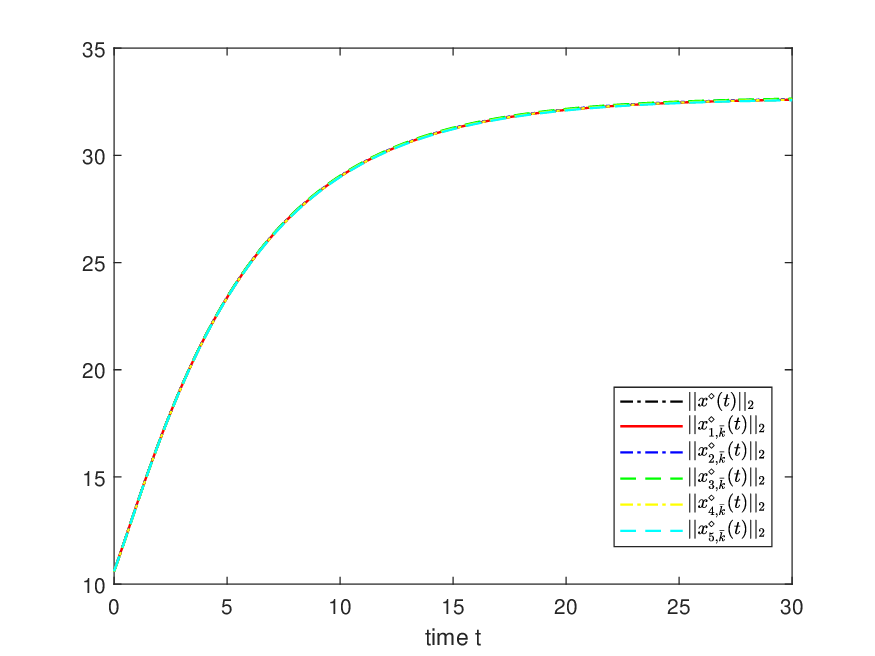}\\
  {\footnotesize Fig. 9. The trajectories of $\|x_{i,\bar{k}}^{\diamond}(t)\|_2$.}
  \end{center}
\end{figure}
\begin{figure}
\begin{center}
  \includegraphics[width=8cm,height=6.3cm]{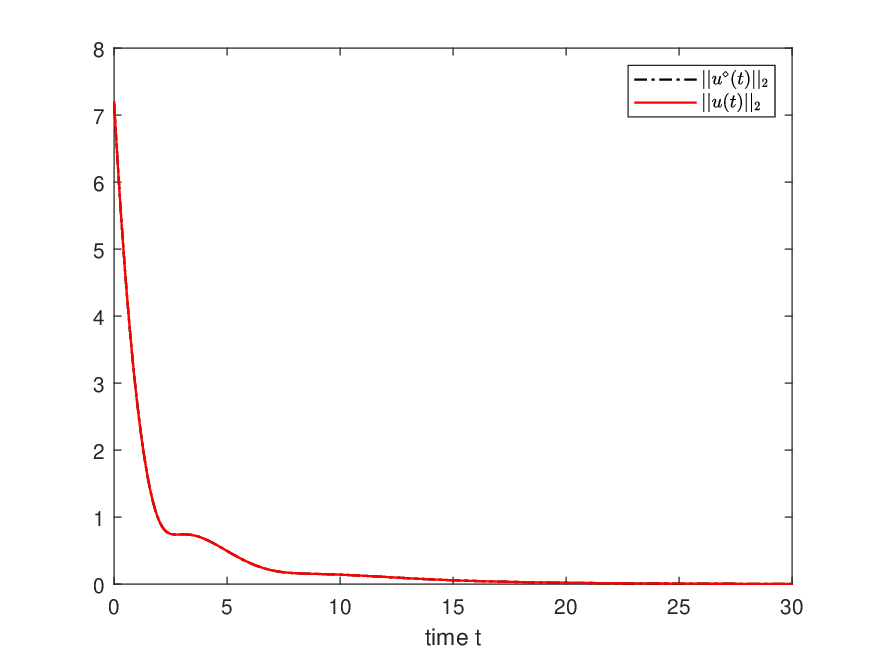}\\
  {\footnotesize Fig. 10. The trajectories of $\|u(t)\|_2$.}
  \end{center}
\end{figure}
\begin{figure}
\begin{center}
  \includegraphics[width=8cm,height=6.3cm]{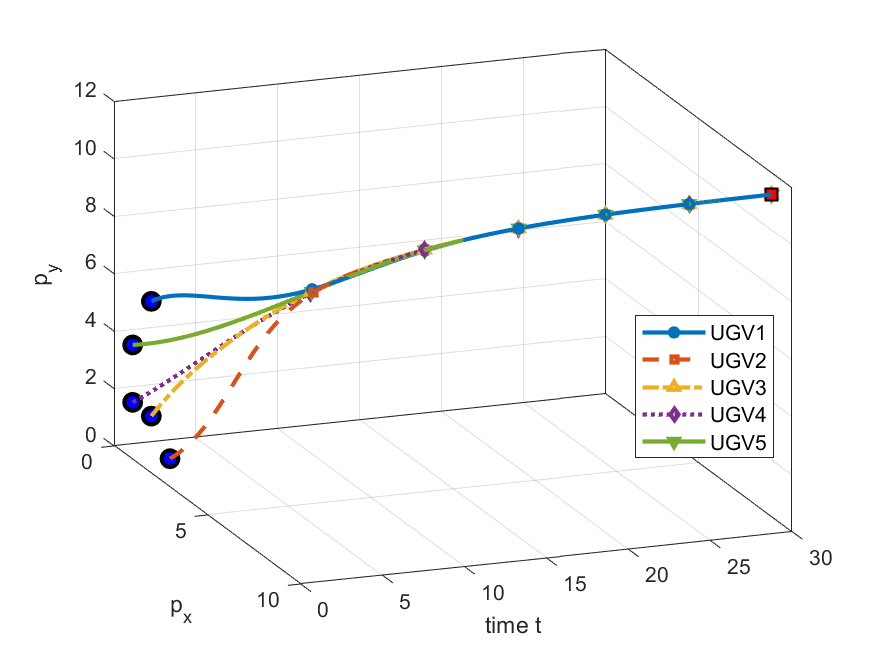}\\
  {\footnotesize Fig. 11. The trajectories of the position states $p_i(t)$ for all UGVs.}
  \end{center}
\end{figure}
\begin{figure}
\begin{center}
  \includegraphics[width=8cm,height=6.3cm]{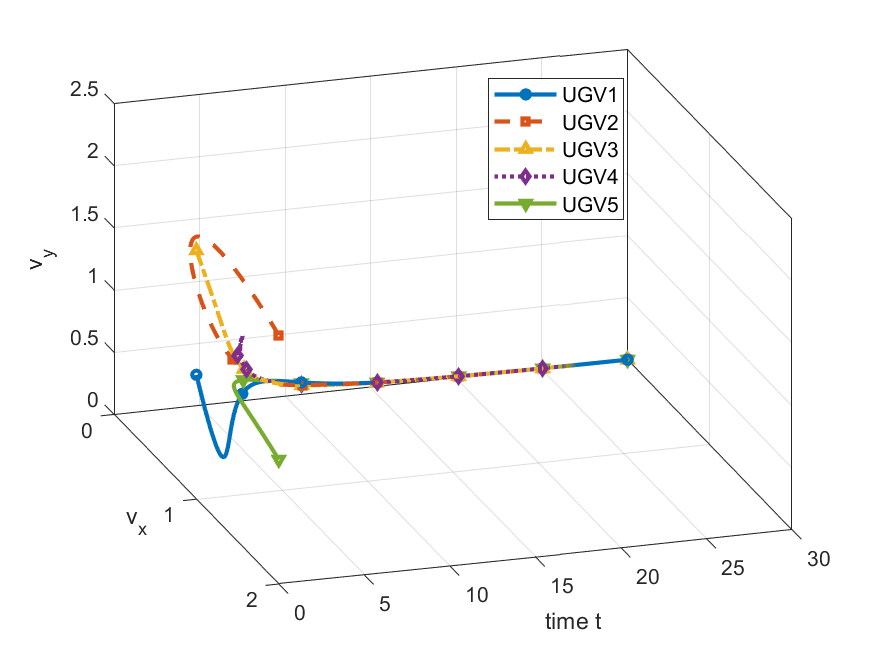}\\
  {\footnotesize Fig. 12. The trajectories of the velocity states $v_i(t)$ for all UGVs.}
  \end{center}
\end{figure}

\section{Conclusions}

This paper studied the LQ optimal control problem for continuous-time systems with terminal state constraints. A distributed algorithm was proposed to solve the optimal controller where only partial information is accessible to the agent, consisting of the distributed iteration of the Riccati equation, the optimal Lagrange multiplier, and the optimal state. As a significant application, the proposed distributed approach is extended to address the optimal consensus control problem. Numerical examples demonstrated the effectiveness of the proposed algorithm, where the performance index under the proposed distributed controller is smaller than that under the commonly used consensus control.

\section*{Appendix A \\ Proof of Lemma 2}\label{secA1}

We will prove the result by using mathematical
induction. For $n=1$, it is obvious from (\ref{13})-(\ref{14}) that
\begin{eqnarray}
Z_{i,k}^{1}(t)&=&Z_{i,k-1}^{1}(t)+\alpha_k[A_{i}-M_iM_i'P_{i,k-1}^{0}(t)\nonumber\\
&&-Z_{i,k-1}^{1}(t)]+\frac{1}{\gamma}\sum_{j\in N_i}[Z_{j,k}^{1}(t)-Z_{i,k-1}^{1}(t)],\label{16}\\
V_{i,k}^{1}(t)&=&V_{i,k-1}^{1}(t)+\alpha_k[Q_{i}+P_{i,k-1}^{0}(t)M_iM_i'P_{i,k-1}^{0}(t)
\nonumber\\
&&-V_{i,k-1}^{1}(t)] +\frac{1}{\gamma}\sum_{j\in N_i}[V_{j,k}^{1}(t)-V_{i,k-1}^{1}(t)].\label{17}
\end{eqnarray}
It is aimed to show
\begin{eqnarray}
\displaystyle\lim_{k\rightarrow \infty}\|Z_{i,k}^{1}(t)-Z^{1}(t)\|&=&0,\label{60}\\
\displaystyle\lim_{k\rightarrow \infty}\|V_{i,k}^{1}(t)-V^{1}(t)\|&=&0.\label{61}
\end{eqnarray}
which is divided into two steps. The first step is to prove (\ref{16})-(\ref{17}) achieves consensus, that is,
\begin{eqnarray}
\displaystyle\lim_{k\rightarrow \infty}\|Z_{i,k}^{1}(t)-Z_{j,k}^{1}(t)\|&=&0,\label{50}\\
\displaystyle\lim_{k\rightarrow \infty}\|V_{i,k}^{1}(t)-V_{j,k}^{1}(t)\|&=&0.\label{51}
\end{eqnarray}
The second step is to prove (\ref{16})-(\ref{17}) achieves global convergence, that is,
\begin{eqnarray}
\displaystyle\lim_{k\rightarrow \infty}\|\frac{1}{N}\sum_{i=1}^{N}Z_{i,k}^{1}(t)-Z^{1}(t)\|&=&0,\label{52}\\
\displaystyle\lim_{k\rightarrow \infty}\|\frac{1}{N}\sum_{i=1}^{N}V_{i,k}^{1}(t)-V^{1}(t)\|&=&0.\label{53}
\end{eqnarray}
Firstly, we prove that the consensus property by denoting
\begin{eqnarray*}
\mathcal{Z}_{k}^{1}(t)&=&\textbf{col}\{Z_{1,k}^{1}(t),Z_{2,k}^{1}(t), \cdots , Z_{N,k}^{1}(t)\},\\
\mathcal{V}_{k}^{1}(t)&=&\textbf{col}\{V_{1,k}^{1}(t),V_{2,k}^{1}(t), \cdots ,
V_{N,k}^{1}(t) \},\\
\Pi_{k}^{1}(t)&=&\textbf{col}\{Y_{1,k}^{1}(t),Y_{2,k}^{1}(t), \cdots ,
Y_{N,k}^{1}(t))\}, \\
\Theta_{k}^{1}(t)&=&\textbf{col}\{L_{1,k}^{1}(t),L_{2,k}^{1}(t), \cdots ,
L_{N,k}^{1}(t)\},\\
Y_{i,k}^{1}(t)&=&A_{i}-M_iM_i'P_{i,k}^{0}(t)-Z_{i,k}^{1}(t), \\ L_{i,k}^{1}(t)&=&Q_{i}+P_{i,k}^{0}(t)M_iM_i'P_{i,k}^{0}(t)-V_{i,k}^{1}(t),
\end{eqnarray*}
it follows from (\ref{16})-(\ref{17}) that
\begin{eqnarray}
\mathcal{Z}_{k}^{1}(t)&=&(\mathcal{A}\otimes I_n)\mathcal{Z}_{k-1}^{1}(t) +\alpha_k\Pi_{k-1}^{1}(t),\nonumber\\
\mathcal{V}_{k}^{1}(t)&=&(\mathcal{A}\otimes I_n)\mathcal{V}_{k-1}^{1}(t) +\alpha_k\Theta_{k-1}^{1}(t),\nonumber
\end{eqnarray}
where $\mathcal{A}=I_N-\frac{1}{\gamma}\mathcal{L}$, $\mathcal{L}$ is the Laplacian matrix \cite{31}. It is obvious that $\mathcal{Z}_{k}^{1}(t)$ and $\mathcal{V}_{k}^{1}(t)$ are bounded. In fact, by defining $\bar{\mathcal{Z}}_{k}^{1}(t)=(U'\otimes I_n)\mathcal{Z}_{k}^{1}(t)$ and $\bar{\mathcal{V}}_{k}^{1}(t)=(U'\otimes I_n)\mathcal{V}_{k}^{1}(t)$, where $U$ is orthogonal matrix such that $U'\mathcal{L}U=diag\{\lambda_1,\cdots,\lambda_N\}$ and $0=\lambda_1\leq\lambda_2\leq\cdots\leq\lambda_N$ are the eigenvalues of $\mathcal{L}$. According to (\ref{16})-(\ref{17}), we have that
\begin{eqnarray}
\bar{\mathcal{Z}}_{k}^{1}(t)&=&[(I-\alpha_{k}I-\frac{1}{\gamma}U'\mathcal{L}U)\otimes I_n]\bar{\mathcal{Z}}_{k-1}^{1}(t)\nonumber\\
&&
-\alpha_{k}(U'\otimes I_n)\Pi_{k-1}^{1}(t),\nonumber\\
\bar{\mathcal{V}}_{k}^{1}(t)&=&[(I-\alpha_{k}I-\frac{1}{\gamma}U'\mathcal{L}U)\otimes I_n]\bar{\mathcal{V}}_{k-1}^{1}(t)\nonumber\\
&&-\alpha_{k}(U'\otimes I_n)\Theta_{k-1}^{1}(t).\nonumber
\end{eqnarray}
Together with the properties of Laplacian matrices \cite{31}, it is easy to get $\bar{\mathcal{Z}}_{k}^{1}(t)$ and $\bar{V}_{k}^{1}(t)$ are bounded. This implies the boundedness of $\mathcal{Z}_{k}^{1}(t)$ and $\mathcal{V}_{k}^{1}(t)$.

Furthemore, let $\mathcal{B}=\frac{1}{N}1_N1_N'$, $\delta_{k}^{1}(t)=[(I_N-\mathcal{B})\otimes I_n]\mathcal{Z}^{1}_{k}(t)$, and $\eta_{k}^{1}(t)=[(I_N-\mathcal{B})\otimes I_n]\mathcal{V}_{k}^{1}(t)$, we have that
\begin{eqnarray}
\delta_{k}^{1}(t)&=&[(\mathcal{A}-\mathcal{B})\otimes I_n]\delta_{k-1}^{1}(t)\nonumber\\
&&+\alpha_k[(I_N-\mathcal{B})\otimes I_n]\Pi_{k-1}^{1}(t),  \label{49}\\
\eta_{k}^{1}(t)&=&[(\mathcal{A}-\mathcal{B})\otimes I_n]\eta_{k-1}^{1}(t)\nonumber\\
&&
+\alpha_k[(I_N-\mathcal{B})\otimes I_n]\Theta_{k-1}^{1}(t),  \label{20}
\end{eqnarray}
where the fact $\mathcal{A}\mathcal{B} = \mathcal{B}\mathcal{A} = \mathcal{B}^{2} = \mathcal{B}$ has been used in the derivations. By applying iteration calculation to (\ref{49}) and (\ref{20}), it yields that
\begin{eqnarray}
\delta_{k}^{1}(t)&=&[(\mathcal{A}-\mathcal{B})\otimes I_n]^{k}\delta_{0}^{1}(t)+\sum_{\tau=1}^{k}\alpha_\tau[(\mathcal{A}-\mathcal{B})\otimes I_n]^{k-\tau}\nonumber\\
&&\times[(I_N-\mathcal{B})\otimes I_n]\Pi_{\tau-1}^{1}(t),  \nonumber\\
\eta_{k}^{1}(t)&=&[(\mathcal{A}-\mathcal{B})\otimes I_n]^{k}\eta_{0}^{1}(t)+\sum_{\tau=1}^{k}\alpha_\tau[(\mathcal{A}-\mathcal{B})\otimes I_n]^{k-\tau}\nonumber\\
&&\times[(I_N-\mathcal{B})\otimes I_n]\Theta_{\tau-1}^{1}(t). \nonumber
\end{eqnarray}
This further gives that
\begin{eqnarray*}
\big\|\delta_{k}^{1}(t)\big\|&=&\big\|[(\mathcal{A}-\mathcal{B})\otimes I_n]^{k}\delta_{0}^{1}(t)+\sum_{\tau=1}^{k}
\alpha_\tau[(\mathcal{A}-\mathcal{B})\nonumber\\
&&\otimes I_n]^{k-\tau}(I_N-\mathcal{B})\otimes I_n]\Pi_{\tau-1}^{1}(t)\big\|\nonumber\\
&\leq&
c\rho^k\|\delta_{0}^{1}(t)\|+\sum_{\tau=1}^{k}\alpha_\tau\|[(\mathcal{A}-\mathcal{B})\otimes I_n]^{k-\tau}\|\nonumber\\
&&\times\|[(I_N-\mathcal{B})\otimes I_n]\|\|\Pi_{\tau-1}^{1}(t)\|,  \nonumber\\
\|\eta_{k}^{1}(t)\|&=&\|[(\mathcal{A}-\mathcal{B})\otimes I_n]^{k}\eta_{0}^{1}(t)+\sum_{\tau=1}^{k}\alpha_\tau[(\mathcal{A}-\mathcal{B})\nonumber\\
&&\otimes I_n]^{k-\tau}[(I_N-\mathcal{B})\otimes I_n]\Theta_{\tau-1}^{1}(t)\|\nonumber\\
&\leq&
c\rho^k\|\eta_{0}^{1}(t)\|+\sum_{\tau=1}^{k}\alpha_\tau\|[(\mathcal{A}-\mathcal{B})\otimes I_n]^{k-\tau}\|\nonumber\\
&&\times\|[(I_N-\mathcal{B})\otimes I_n]\|\|\Theta_{\tau-1}^{1}(t)\|. \nonumber
\end{eqnarray*}
Since the graph is connected, there exists a constant $\gamma>0$ such that the following inequality holds
\begin{eqnarray*}
\|[(\mathcal{A}-\mathcal{B})\otimes I_n]^{k}\leq c\rho^k,
\end{eqnarray*}
where $c > 0$ and $\rho \in (0, 1)$. According to $\|[(I_N-\mathcal{B})\otimes I_n]\| < \infty$,
$\|\delta_{0}^{1}(t)\| < \infty$, $\|\eta_{0}^{1}(t)\| < \infty$, $\|\Pi_{\tau}^{1}(t)\| < \infty$,
$\|\Theta_{\tau}^{1}(t)\| < \infty$, we have that
\begin{eqnarray}
&&\lim_{k\to \infty} c\rho^k\|\delta_{0}^{1}\|=0,\label{40}\\
&&\lim_{k\to \infty} c\rho^k\|\eta_{0}^{1}\|=0,\label{41}\\
&&\lim_{k\to\infty}\sum_{\tau=1}^{k}\alpha_\tau\|[(\mathcal{A}-\mathcal{B})\otimes I_n]^{k-\tau}\|\nonumber \\
&&\quad\quad\quad\quad\quad\times\|[(I_N-\mathcal{B})\otimes I_n]\|\|\Pi_{\tau-1}^{1}(t)\|=0,\label{42}\\
&&\lim_{k\to \infty}\sum_{\tau=1}^{k}\alpha_\tau\|[(\mathcal{A}-\mathcal{B})\otimes I_n]^{k-\tau}\|\nonumber \\
&&\quad\quad\quad\quad\quad\times\|[(I_N-\mathcal{B})\otimes I_n]\|\|\Theta_{\tau-1}^{1}(t)\|=0.\label{43}
\end{eqnarray}

Thus, we obtain that
\begin{eqnarray*}
\lim_{k\to \infty}\|\delta_{k}^{1}(t)\| &=& 0,\\
\lim_{k\to \infty}\|\eta_{k}^{1}(t)\|&= &0,
\end{eqnarray*}
that is, (\ref{50})-(\ref{51}) hold.

Secondly, we prove that (\ref{52})-(\ref{53}) hold. To do this,
we define $\tilde{Z}_{k}^{1}(t)=\frac{1}{N}\sum_{i=1}^{N}Z_{i,k}^{1}(t)$ and $\tilde{V}_{k}^{1}(t)=\frac{1}{N}\sum_{i=1}^{N}V_{i,k}^{1}(t)$, and derive from (\ref{16})-(\ref{17}) that
\begin{eqnarray}
 \tilde{Z}_{k}^{1}(t)
 &=&\tilde{Z}_{k-1}^{1}(t)+\frac{1}{N}\alpha_k\sum_{i=1}^{N}[A_{i}-M_iM_i'P_{i,k-1}^{0}(t)
 \nonumber\\
&&-Z_{i,k-1}^{1}(t)],\nonumber\\
  \tilde{V}_{k}^{1}(t)
 &=&\tilde{V}_{k-1}^{1}(t)+\frac{1}{N}\alpha_k\sum_{i=1}^{N}[Q_{i}+P_{i,k-1}^{0}(t)M_iM_i'
\nonumber\\
&&\times P_{i,k-1}^{0}(t)-V_{i,k-1}^{1}(t)].\nonumber
\end{eqnarray}

By defining  $Z_{k}^{1}(t)$, $V_{k}^{1}(t)$ by
\begin{eqnarray}
Z_{k}^{n}(t)&=&Z_{k-1}^{n}(t)+\alpha_k[A-BR^{-1}B'P^{n-1}(t)\nonumber \\
&&-Z_{k-1}^{n}(t)],\label{10}\\
V_{k}^{n}(t)&=&V_{k-1}^{n}(t)+\alpha_k[Q+P^{n-1}(t)BR^{-1}B'\nonumber \\
&&
\times P^{n-1}(t)-V_{k-1}^{n}(t)],\label{11}
\end{eqnarray}
and letting $\Delta Z_k^{1}(t) = \tilde{Z}_{k}^{1}(t) - Z_{k}^{1}(t) $, $\Delta V_k^{1}(t) = \tilde{V}_{k}^{1}(t) - V_{k}^{1}(t)$, we have that
\begin{eqnarray}
\Delta{Z}_{k}^{1}(t)
 &=&(1-\alpha_k)\Delta{Z}_{k-1}^{1}(t),\nonumber\\
\Delta{V}_{k}^{1}(t)
 &=&(1-\alpha_k)\Delta{V}_{k-1}^{1}(t).\nonumber
\end{eqnarray}

By applying the facts that $0<1-\alpha_k<1$, and $\lim_{k\to \infty} \alpha_k = 0$, it yields that
\begin{eqnarray}
\lim_{k\to \infty} \|\Delta{Z}_{k}^{1}(t)\|&=&0,\label{45}\\
\lim_{k\to \infty} \|\Delta{V}_{k}^{1}(t)\|&=&0.\label{48}
\end{eqnarray}

Combining with (\ref{50})-(\ref{51}), it follows that
\begin{eqnarray*}
\lim_{k\to \infty} \|Z_{i,k}^{1}(t) - Z_{k}^{1}(t) \|&=&0,\\
\lim_{k\to \infty} \|V_{i,k}^{1}(t) - V_{k}^{1}(t) \|&=&0.
\end{eqnarray*}

Moreover, together with the fact that $\lim_{k\rightarrow \infty}\|Z_{k}^{1}(t)-Z^{1}(t)\| = 0$ and $\lim_{k\rightarrow \infty}\|V_{k}^{1}(t)-V^{1}(t)\| = 0$ obtained by using Theorem 3.1.1 in \cite{36},  we obtain (\ref{60})-(\ref{61}) directly.

To complete the induction proof, assume that for any $m\leq n$, (\ref{19})-(\ref{21}) holds, and $P_{i,k}^{m}(t)$ obtained by (\ref{31}) is convergent. By taking similar proof for (\ref{50})-(\ref{51}) and (\ref{52})-(\ref{53}), it can be obtained that (\ref{19})-(\ref{21}) holds when $n=m+1$. To avoid repetition, we omit the details of the proof here. \hfill $\blacksquare$


\begin{thebibliography}{0}

\bibitem{1}
I. Munteanu, N. A. Cutululis, A. I. Bratcu, and E. Ceang$\breve{a}$,
{``Optimization of variable speed wind power systems based on a LQG approach,''}
{\em Control engineering practice,} {vol. 13, no. 4, pp. 903-912, 2005.}

\bibitem{2}
P. Benigno, and M. Woodford,
{``Optimal monetary and fiscal policy: A linear-quadratic approach,''}
{\em NBER macroeconomics annual,} {vol. 18, pp. 271-333, 2003.}


\bibitem{3}
N. E. Kahveci, P. A. Ioannou, and M. D. Mirmirani,
{``Adaptive LQ control with anti-windup augmentation to optimize UAV performance in autonomous soaring applications,''}
{\em IEEE Transactions on Control Systems Technology,} {vol. 16, no. 4, pp. 691-707, 2008.}


\bibitem{4}
Q. Sun, J. Xu, and H. Zhang.
{``Guidance for hypersonic reentry using nonlinear model predictive control and radau pseudospectral method,''}
{\em Optimal Control Applications and Methods,} {DOI: 10.1002/oca.3267.}

\bibitem{5}
R. Kalman,
{``Contributions to the theory of optimal control,''}
{\em Bol. Soc. Mat. Mexicana,} {vol. 5, no. 63, pp. 102-119, 1960.}

\bibitem{6}
H. Zhang, and J. Xu,
{``Optimal control with irregular performance,''}
{\em Science China Information Sciences,} {vol. 62, no. 9, pp. 1-14, 2019.}

\bibitem{7}
J. Zhao, and R. Zhou,
{``Reentry trajectory optimization for hypersonic vehicle satisfying complex constraints,''}
{\em Chinese Journal of Aeronautics,} {vol. 26, no. 6, pp. 1544-1553, 2013.}

\bibitem{8}
H. Michalska, and D. Q. Mayne,
{``Robust receding horizon control of constrained nonlinear systems,''}
{\em IEEE Transactions on Automatic Control,} {vol. 38, no. 11, pp. 1623-1633, 1993.}

\bibitem{9}
S. R. Vadali and R. Sharma,
{``Optimal finite-time feedback controllers for nonlinear systems with terminal constraints,''}
{\em Journal of Guidance, Control, and Dynamics,} {vol. 29, no. 4, pp. 921-928, 2006.}

\bibitem{10}
X. Bi, J. Sun, and J. Xiong.
{``Optimal control for controllable stochastic linear systems,''}
{\em ESAIM: Control, Optimisation and Calculus of Variations,} {vol. 26, no. 98, 2020.}


\bibitem{11}
S. Ji, and X. Zhou,
{``A maximum principle for stochastic optimal control with terminal state constraints, and its applications,''}
{\em Communications in Information $\&$ Systems,} {vol. 6, no. 4, pp. 321-338, 2006.}


\bibitem{12}
J. Sun,
{``Linear quadratic optimal control problems with fixed terminal states and integral quadratic constraints,''}
{\em  Applied Mathematics $\&$ Optimization,} {vol. 83, no. 1, pp. 251-276, 2021.}


\bibitem{13}
J. Liu, J. Xu, H. Zhang, and M. Fu,
{``Stochastic LQ optimal control with initial and terminal constraints,''}
{\em IEEE Transactions on Automatic Control,} {vol. 69, no. 9, pp. 6261-6268, 2024.}

\bibitem{14}
M. Tubaishat, and S. Madria,
{``Sensor networks: an overview,''}
{\em IEEE Potentials,} {vol. 22, no. 2, pp. 20-23, 2003.}



\bibitem{15}
X. Li, W. Shu, M. Li, H. -Y. Huang, P. -E. Luo, and M. -Y. Wu,
{``Performance evaluation of vehicle-based mobile sensor networks for traffic monitoring,''}
{\em IEEE Transactions on Vehicular Technology,} {vol. 58, no. 4, pp. 1647-1653, 2009.}

\bibitem{16}
V. C. Gungor, B. Lu, and G. P. Hancke,
{``Opportunities and challenges of wireless sensor networks in smart grid,''}
{\em IEEE Transactions on Industrial Electronics,} {vol. 57, no. 10, pp. 3557-3564, 2010.}

\bibitem{17}
L. M. Oliveira, and J. J. Rodrigues,
{``Wireless Sensor Networks: A survey on environmental monitoring,''}
{\em Journal Of Communications,} {vol. 6, no. 2, pp. 143-151, 2011.}

\bibitem{18}
G. Zhao,
{``Wireless sensor networks for industrial process monitoring and control: A survey,''}
{\em Netw. Protocols Algorithms,} {vol. 3, no. 1, pp. 46-63, 2011.}

\bibitem{19}
T. Ba\c{s}ar,
{``Two-criteria LQG decision problems with one-step delay observation sharing pattern,''}
{\em Information and Control,} {vol. 38, no. 1, pp. 21-50, 1978.}

\bibitem{20}
X. Liang, J. Xu, H. Wang, and H. Zhang,
{``Decentralized output-feedback control with asymmetric one-step delayed information,''}
{\em IEEE Transactions on Automatic Control,} {vol. 68, no. 12, pp. 7871-7878, 2023.}

\bibitem{21}
A. Rantzer,
{``Linear quadratic team theory revisited,''}
{\em Procedure 2006 American Control Conference, Minneapolis,} {2006, pp. 1637-1641.}


\bibitem{22}
J. Bismut,
{``An example of interaction between information and control: The Transparency of a game,''}
{\em IEEE Transactions on Automatic Control,} {vol. 18, no. 5, pp. 518-522, 1973.}

\bibitem{23}
J. Wu, and S. Lall,
{``A dynamic programming algorithm for decentralized Markov decision processes with a broadcast structure,''}
{\em Procedure 49th IEEE Conference on Decision and Control,} {2010, pp. 6143-6148.}

\bibitem{24}
A. Mahajan, A. Nayyar, and D. Teneketzis,
{``Identifying tractable decentralized control problems on the basis of information structure,''}
{\em Procedure 46th Annual Allerton Conference on Communication, Control, and Computing,} {Monticello, IL, USA, 2008, pp. 1440-1449.}

\bibitem{25}
M. A. Rami, X. Chen, and X. Y. Zhou,
{``Coordination of groups of mobile autonomous agents using nearest neighbor rules,''}
{\em IEEE Transactions on Automatic Control,} {vol. 48, no. 6, pp. 988-1001, 2003.}

\bibitem{26}
G. Gu, L. Marinovici, and F. L. Lewis,
{``Consensusability of discrete-time dynamic multiagent systems,''}
{\em IEEE Transactions on Automatic Control,} {vol. 57, no. 8, pp. 2085-2089, 2012.}

\bibitem{27}
W. Wang, F. Zhang, and C. Han,
{``Distributed LQR control for discrete-time homogeneous systems,''}
{\em International Journal of Systems Science,} {vol. 47, no. 15, pp. 3678-3687, 2016.}

\bibitem{28}
Z. Zhang, W. Yan, and H. Li,
{``Distributed optimal control for linear multiagent systems on general digraphs,''}
{\em IEEE Transactions on Automatic Control,} {vol. 66, no. 1, pp. 322-328, 2021.}

\bibitem{29}
F. Borrelli, and T. Keviczky,
{``Distributed LQR design for identical dynamically decoupled systems,''}
{\em IEEE Transactions on Automatic Control,} {vol. 53, no. 8, pp. 1901-1912, 2008.}

\bibitem{30}
S. Yao, S. Xie, and T. Li,
{``Online distributed optimization algorithm with dynamic regret analysis under unbalanced graphs,''}
{\em Automatica,} {vol. 174, no. 112116, 2025.}

\bibitem{31}
S. Battilotti, F. Cacace, and M. d'Angelo,
{``Distributed optimal control of discrete-time linear systems over networks,''}
{\em IEEE Transactions on Control of Network Systems,} {vol. 11, no. 2, pp. 671-682, 2024.}

\bibitem{32}
J. Li, G. Deng, C. Luo, Q. Lin, Q. Yan, and Z. Ming,
{``A hybrid path planning method in unmanned air/ground vehicle (UAV/UGV) cooperative systems,''}
{\em IEEE Transactions on Vehicular Technology,} {vol. 65, no. 12, pp. 9585-9596, 2016.}

\bibitem{33}
J. Xu, J. Liu, Z. Zhang, and W. Wang,
{``Q-Learning for linear quadratic optimal control with terminal state constraint,''}
{\em Optimization and Control,} {vol. 29, pp. 134-140, 2024.}

\bibitem{34}
J. Yong, and X. Zhou,
{``Stochastic controls: Hamiltonian systems and HJB equations,''}
{\em Springer Science $\&$ Business Media,} {vol. 43, 2012.}

\bibitem{35}
L. An, and G. -H. Yang,
{``Distributed optimal coordination for heterogeneous linear multiagent systems,''}
{\em IEEE Transactions on Automatic Control,} {vol. 67, no. 12, pp. 6850-6857, 2022.}

\bibitem{36}
H. Chen,
{``Stochastic approximation and its applications,''}
{ Norwell, MA: Kluwer, 2002.}

\bibitem{37}
C.-Q. Ma, and J.-F. Zhang,
{`` Necessary and sufficient conditions for consensusability of linear multi-agent systems,''}
{\em IEEE Transactions on Automatic Control,} {vol. 55, no. 5, pp. 1263-1268, 2010.}


\bibitem{38}
R. O. -Saber, and R. Murray,
{``Consensus problems in networks of agents with switching topology and time-delays,''}
{\em IEEE Transactions on Automatic Control,} {vol. 49, no. 9, pp. 1520-1533, 2004.}

\bibitem{39}
W. Ren, and R. Beard,
{``Consensus seeking in multiagent systems under dynamically changing interaction topologies,''}
{\em IEEE Transactions on Automatic Control,} {vol. 50, no. 5, pp. 655-661, 2005.}


\end{thebibliography}
 \end{document}